%% file: main-1.tex
\documentclass[final]{siamart171218}
\input{ex_shared}

\ifpdf
\hypersetup{
  pdftitle={\TheTitle},
  pdfauthor={\TheAuthors}
}
\fi

\begin{document}
\maketitle

\begin{abstract}
We present a general technique for the analysis of first-order methods. The technique relies on the construction of a duality gap for an appropriate approximation of the objective function, where the function approximation improves as the algorithm converges. We show that in continuous time enforcement of an invariant that this approximate duality gap decreases at a certain rate exactly recovers a wide range of first-order continuous-time methods. We characterize the discretization errors incurred by different discretization methods, and show how iteration-complexity-optimal methods for various classes of problems cancel out the discretization error. The techniques are illustrated on various classes of problems -- including convex minimization for Lipschitz-continuous objectives, smooth convex minimization, composite minimization, smooth and strongly convex minimization, solving variational inequalities with monotone operators, and convex-concave saddle-point optimization -- and naturally extend to other settings.   
\end{abstract}



\section{Introduction}\label{sec:intro}

First-order optimization methods have recently gained high popularity due to their applicability to large-scale problem instances arising from modern datasets, their relatively low computational complexity, and their potential for parallelizing computation \cite{sra2012optimization}. Moreover, such methods have also been successfully applied in discrete optimization leading to faster  numerical methods~\cite{ST04,KOSZ13}, graph algorithms~\cite{KLOS2014,Sherman2013,LRS2013},  and submodular optimization methods~\cite{Ene}.

Most first-order optimization methods can be obtained from the discretization of continuous-time dynamical systems that converge to optimal solutions. In the case of mirror descent, the continuous-time view was the original motivation for the algorithm~\cite{nemirovskii1983problem}, while more recent work has focused on deducing continuous-time interpretations of accelerated methods~\cite{wibisono2016variational,wilson2016lyapunov,krichene2015accelerated,SuBC16,Scieur2017}. 

Motivated by these works, we focus on providing a unifying theory of first-order methods as discretizations of continuous-time dynamical systems. We term this general framework  
the \emph{Approximate Duality Gap Technique (\adgt)}. In addition to providing an intuitive and unified convergence analysis of various first-order methods that is often only a few lines long, \adgt is also valuable in developing new first-order methods with tight convergence bounds%
~\cite{AXGD,cohen2018acceleration}, {in clarifying interactions between the acceleration and noise~\cite{cohen2018acceleration}, and in obtaining width-independent\footnote{{Width-independent algorithms enjoy poly-logarithmic dependence of their convergence times on the constraints matrix width -- namely, the ratio between the constraint matrix maximum and minimum non-zero elements. By contrast, standard first-order methods incur (at best) linear dependence on the matrix width, which is not even considered to be polynomial-time convergence~\cite{nesterov2005smooth}.}} algorithms for problems with positive linear constraints~\cite{LP-jelena-lorenzo,diakonikolas2018fairpc}. Further, we have extended \adgt to the setting of block coordinate descent methods~\cite{diakonikolas2018alternating}.} 

Unlike traditional approaches that start from an algorithm description and then use arguments such as 
Lyapunov stability criteria to prove convergence bounds~\cite{nemirovskii1983problem,wibisono2016variational,wilson2016lyapunov,krichene2015accelerated,SuBC16}, 
our approach takes the opposite direction: \emph{continuous-time methods are obtained from the analysis}, using purely optimization-motivated arguments.  

{In particular, \adgt can be summarized as follows. Given a convex optimization problem $\min_{\vx \in X}f(\vx)$, to show that a method converges to a minimizer $\vx^* \in \argmin_{\vx \in X} f(\vx)$ at rate $1/\a^{(t)}$ (e.g., $\a^{(t)} = t$ or $\a^{(t)}=t^2$), we need to show that $f(\vx^{(t)}) - f(\vx^*) \leq \nicefrac{Q}{\a^{(t)}}$, where $\vx^{(t)}\in X$ is the solution produced by the method at time $t$ and $Q \in \mathbb{R_+}$ is some bounded quantity that is independent of time $t$. In general, keeping track of the true optimality gap $f(\vx^{(t)}) - f(\vx^*)$ is challenging, as the minimum function value $f(\vx^*)$ is typically not known to the method. Instead, the main idea of \adgt is to create an estimate of the optimality gap $G^{(t)}$ that can be easily tracked and controlled and ensure that $\a^{(t)}G^{(t)}$ is a non-increasing function of time. The estimate corresponds to the difference between an upper bound on $f(\vx^{(t)})$, $U^{(t)}\geq f(\vx^{(t)}),$ and a lower bound on $f(\vx^*),$ $L^{(t)}\leq f(\vx^*),$ so that $f(\vx^{(t)}) - f(\vx^*)\leq U^{(t)}-L^{(t)} = G^{(t)}$. Since \adgt ensures that $\a^{(t)}G^{(t)}$ is a non-increasing function of time, it follows that $f(\vx^{(t)}) - f(\vx^*)\leq G^{(t)} \leq \nicefrac{\alpha^{(0)}G^{(0)}}{\alpha^{(t)}}$, which is precisely what we want to show, as long as we ensure that $\alpha^{(0)}G^{(0)}$ is bounded.}

To illustrate the power and generality of the technique, we show how to obtain and analyze several well-known first-order methods, such as gradient descent, dual averaging~\cite{nesterov2009primal}, mirror-prox/extra-gradient method~\cite{extragradient-descent,Nesterov2007dual-extrapolation,Mirror-Prox-Nemirovski}, accelerated methods~\cite{Nesterov1983,nesterov2005smooth}, composite minimization methods \cite{duchi2010composite,nesterov2015universal}, and Frank-Wolfe methods~\cite{nesterov2015cgm}. The same ideas naturally extend to other classes of convex optimization problems and their corresponding optimal first-order methods. Here, ``optimal'' is in the sense that the methods yield worst-case iteration complexity bounds for which there is a matching lower bound (i.e., ``optimal'' is in terms of iteration complexity).
\subsection{Related Work}\label{sec:related-work}

There exists a large body of research in optimization and  first-order methods in particular, and, while we cannot provide a thorough literature review, we refer the reader to recent books \cite{sra2012optimization,Bubeck2015,ben2001lectures,nesterov2013introductory}. 

Multiple approaches to unifying analyses of first-order methods have been developed, with a particular focus on explaining the acceleration phenomenon. Tseng gives a formal framework that unifies all the different instantiations of accelerated gradient methods~\cite{tseng2008}. More recently, Allen-Zhu and Orecchia~\cite{AO-survey-nesterov} interpret acceleration as coupling of mirror descent and gradient descent steps. Bubeck~\etal~provide an elegant geometric interpretation of the Euclidean instantiation of Nesterov's method~\cite{Bubeck2015}. Drusvyatskiy~\etal~\cite{drusvyatskiy2016optimal} interpret the geometric descent of Bubeck~\etal~\cite{Bubeck2015} as a sequence minimizing quadratic lower-models of the objective function and obtain limited-memory extensions with improved performance. Lin~\etal~\cite{lin2015universal} provide a universal scheme for accelerating non-accelerated first-order methods.

Su~\etal~\cite{SuBC16} and Krichene~\etal~\cite{krichene2015accelerated} interpret Nesterov's accelerated method as a discretization of a certain continuous-time dynamics and analyze it using Lyapunov stability criteria.  Scieur~\etal~\cite{Scieur2017} interpret acceleration as a multi-step integration method from numerical analysis applied to the gradient flow differential equation. 
Wibisono~\etal~\cite{wibisono2016variational} and Wilson~\etal~\cite{wilson2016lyapunov} interpret accelerated methods using Lyapunov stability analysis and drawing ideas from Lagrangian mechanics. 
A recent note of Bansal and Gupta~\cite{bansal2017potential} provides an intuitive, potential-based interpretation of many of the commonly-known convergence proofs for first-order methods.

The two references~\cite{wibisono2016variational,wilson2016lyapunov} are most closely related to our work, as they are motivated by the Lagrangian view from classical mechanics that leads to the description of system dynamics through the principle of stationary action. In a similar spirit, most dynamics described in our work maintain the invariant that $\frac{\d}{\d t}(\a^{(t)}G^{(t)})=0.$  However, unlike our work, which relies on enforcing an invariant that both leads to the algorithms and their convergence analysis, the results from~\cite{wibisono2016variational,wilson2016lyapunov} rely on the use of separate Lyapunov functions to obtain convergence results. It is unclear how these Lyapunov functions relate to the problems' optimality gaps, which makes them harder to generalize, especially in non-standard settings such as, e.g.,~\cite{LP-jelena-lorenzo,diakonikolas2018fairpc,diakonikolas2018alternating}.

The approximate duality gap presented here is closely related to Nesterov's estimate sequence (see e.g.,~\cite{nesterov2013introductory}). In particular, up to the regularization term $\phi(\vx^*)/\a^{(t)}$, our lower bound $L^{(t)}$ is equivalent to Nesterov's estimate sequence, providing a natural interpretation of this powerful and commonly used technique.

\subsection{Notation}\label{sec:notation}
We use non-boldface letters to denote scalars and boldface letters to denote vectors. Superscript index $(\cdot)^{(t)}$ denotes the value of $(\cdot)$ at time $t$. The ``dot'' notation is used to denote the time derivative, i.e., $\dot{x} = \frac{d x}{d t}$. Given a measure $\alpha^{(\tau)}$ defined on $\tau \in [0, t]$, we use the Lebesgue-Stieltjes notation for the integral. In particular, given $\phi^{(\tau)}$ defined on $\tau\in [0, t]$:
\begin{equation*}
\int_{0}^{t} \phi^{(\tau)} \dot{\alpha}^{(\tau)}d\tau = \int_{0}^{t} \phi^{(\tau)} \mathrm{d}\alpha^{(\tau)}. 
\end{equation*}
{
In the discrete-time setting, we will assume that $\a^{(\tau)}$ is an increasing piecewise constant function, with discontinuities occurring only at discrete time points $i \in Z_+$, and such that $\a^{(t)}=0$ for $t <0$. Hence, $\dot{\a}^{(\tau)}$ can be expressed as a train of Dirac Delta functions: $\dot{\a}^{(\tau)} = \sum_{i=0}^\infty a_i\delta(\tau-i),$ where $a_i = \a^{(i+\Delta)} - \a^{(i-\Delta)}$ for $\Delta \in (0, 1)$. This means that $\dot{\a}^{(\tau)}$ samples the function under the integral, so that $\int_{0}^t \phi^{(\tau)}\d\a^{(\tau)} = \sum_{i=0}^{\lfloor t \rfloor}a_i\phi^{(i)}$.} 
We denote $A^{(t)} = \int_{0}^t \mathrm{d}\alpha^{(\tau)}$, so that $\frac{1}{A^{(t)}}\int_{0}^t \mathrm{d}\alpha^{(\tau)} = 1$. {In continuous time, $A^{(t)} = \a^{(t)}-\a^{(0)},$ while in the discrete time $A^{(t)} = \sum_{i=0}^{\lfloor t \rfloor}a_i=\a^{(t)}.$} 
We assume throughout the paper that {$\a^{(0)}>0$ and} $\dot{\alpha}^{(t)}>0, \forall t\geq 0$, and use the following notation for the aggregated negative gradients: 
\begin{equation}\label{eq:def-of-z}
\vz^{(t)}\defeq -\int_{0}^t \nabla f(\vx^{(\tau)})\mathrm{d}\alpha^{(\tau)}.
\end{equation}

For all considered problems, we assume that the feasible region is a closed convex set $X\subseteq \mathbb{R}^n$, for a finite $n$. We assume that there is a (fixed) norm $\|\cdot\|$ associated with $X$ and define its dual norm in a standard way: $\|\vz\|_* = \max_{\vx \in X}\{\innp{\vz, \vx}: \|\vx\|\leq 1\}$, where $\innp{\cdot, \cdot}$ denotes the inner product.
\subsection{Preliminaries}\label{sec:prelims}

We focus on minimizing a continuous and differentiable\footnote{The differentiability assumption is not always necessary and can be relaxed to subdifferentiability in the case of dual averaging/mirror descent methods. Nevertheless, we will assume differentiability throughout the paper for simplicity of exposition.} convex function $f(\cdot)$ defined on a convex set $X \subseteq \mathbb{R}^n$, and we let $\vx^* = \arg\min_{\vx \in X}f(\vx)$ denote the minimizer of $f(\cdot)$ on $X$. The following definitions will be useful in our analysis, and thus we state them here for completeness.  

\begin{definition}\label{def:convexity}
A function $f:X\rightarrow \mathbb{R}$ is convex on $X$, if for all $\vx, \vxh \in X$: $f(\vxh) \geq f(\vx) + \innp{\nabla f(\vx), \vxh - \vx}$.
\end{definition}

\begin{definition}\label{def:smoothness}
A function $f:X\rightarrow \mathbb{R}$ is smooth on $X$ with smoothness parameter $L$ and with respect to a norm $\|\cdot\|$, if for all $\vx, \vxh \in X$: $f(\vxh) \leq f(\vx) + \innp{\nabla f(\vx), \vxh - \vx} + \frac{L}{2}\|\vxh - \vx\|^2$. Equivalently: $\|\nabla f(\vx)-\nabla f(\vxh)\|_*\leq L\|\vx - \vxh\|$.
\end{definition}

\begin{definition}\label{def:strong-convexity}
A function $f:X\rightarrow \mathbb{R}$ is strongly convex on $X$ with strong convexity parameter $\sigma$ and with respect to a norm $\|\cdot\|$, if for all $\vx, \vxh \in X$: $f(\vxh) \geq f(\vx) + \innp{\nabla f(\vx), \vxh - \vx} + \frac{\sigma}{2}\|\vxh - \vx\|^2$. Equivalently: $\|\nabla f(\vx)-\nabla f(\vxh)\|_*\geq \sigma\|\vx - \vxh\|$.
\end{definition}

\begin{definition}\label{def:bregman-divergence}(Bregman Divergence) Bregman divergence of a function $\psi$ is defined as: 
$D_{\psi}(\vx, \vxh) \defeq \psi(\vx) - \psi(\vxh)-\innp{\nabla \psi(\vxh), \vx - \vxh}$. 
\end{definition}

\begin{definition}\label{def:cxv-conj}(Convex Conjugate) Function $\psi^*(\cdot)$ is the convex conjugate of $\psi: X \rightarrow \mathbb{R}$, if $\psi^* (\vz) = \sup_{\vx \in X}\{\innp{\vz, \vx} - \psi(\vx)\}$, $\forall \vz \in \mathbb{R}$. 
\end{definition}
As $X$ is assumed to be closed, $\sup$ in Definition~\ref{def:cxv-conj} can be replaced by $\max$.

We will assume that there is a differentiable strictly convex function $\phi:X\rightarrow \mathbb{R}$, possibly dependent on $t$ (in which case we denote it as $\phi_t$), such that $\max_{\vx \in X} \{\innp{\vz, \vx} - \phi (\vx)\}$ is easily solvable, possibly in a closed form. Notice that this problem defines the convex conjugate of $\phi(\cdot)$, i.e., $\phi^* (\vz) = \max_{\vx \in X} \{\innp{\vz, \vx} - \phi (\vx)\}$. {We will further assume without loss of generality that $\min_{\vx \in X} \phi(\vx) \geq 0$.\footnote{{This assumption can be easily satisfied by taking $\phi(\cdot)$ to be, for example, a Bregman divergence: $\phi(\vx) = D_{\psi}(\vx, \vx^{(0)})$ for some strictly convex and differentiable $\psi$ and fixed $\vx^{(0)} \in X$.}} The role of function $\phi$ will be to regularize the lower bound in the construction of the approximate duality gap.} 
The following standard fact, based on Danskin's Theorem (see, e.g.~\cite{bertsekas1971control,Bertsekas2003}), will be extremely useful in carrying out the analysis of the algorithms in this paper.

\begin{fact}\label{fact:danskin}
Let $\phi: X \to \mathbb{R}$ be a differentiable strongly convex function. Then:
$$
\nabla \phi^*(\vz) = \arg\max_{\vx \in X} \left\{ \innp{\vz, \vx} - \phi(\vx)\right\} = \arg\min_{\vx \in X} \left\{ -\innp{\vz, \vx} + \phi(\vx)\right\}.
$$ 
\end{fact}
In particular, Fact \ref{fact:danskin} implies:
\begin{equation}\label{eq:mirror-map-z-t}
\nabla\phi^*(\vz^{(t)}) = \arg\min_{\vx \in X}\left\{\int_{0}^t \innp{\nabla f(\vx^{(\tau)}), \vx - \vx^{(\tau)}}\mathrm{d}\alpha^{(\tau)} +\phi(\vx)\right\}
\end{equation}

Some other useful properties of Bregman divergence can be found in Appendix~\ref{sec:breg-div-prop}.
\paragraph{Overview of Continuous-Time Operations} 

In continuous time, changes in the variables are described by differential equations. Of particular interest are (weighted) aggregation and averaging. Aggregation of a function $g(x)$ is $\dot{y}^{(t)} = \dot{\alpha}^{(t)}g(x^{(t)})$. 
Observe that, by integrating both sides from $0$ to $t$, this is equivalent to: $y^{(t)} = y^{(0)} + \int_{0}^t g(x^{(\tau)})\mathrm{d}\alpha^{(\tau)}$. Averaging of a function $g(x)$ is $\dot{y}^{(t)} = \dot{\alpha}^{(t)}\frac{g(x^{(t)})-y^{(t)}}{\alpha^{(t)}}$. This can be equivalently written as $\frac{\d}{\d t}(\alpha^{(t)}y^{(t)}) = \dot{\alpha}^{(t)}g(x^{(t)})$, implying $y^{(t)} = \frac{\alpha^{(0)}}{\alpha^{(t)}}y^{(0)} + \frac{1}{\alpha^{(t)}}\int_{0}^t g(\vx^{(\tau)})\mathrm{d}\alpha^{(\tau)}$. 
The following simple proposition will be useful in our analysis.
\begin{proposition}\label{prop:ct-differentiation-inside-min}
$\frac{\d}{\d t} \min_{\vx \in X}\left\{ -\innp{\vz^{(t)}, \vx} + \phi(\vx) \right\} = -\innp{\dot{\vz}^{(t)}, \nabla\phi^*(\vz^{(t)})}$.
\end{proposition}
\begin{proof}
Follows by observing that $\phi^*(\vz^{(t)}) = -  \min_{\vx \in X}\left\{ -\innp{\vz^{(t)}, \vx} + \phi(\vx) \right\}$ and applying the chain rule.
\end{proof}
\section{The Approximate Duality Gap Technique}\label{sec:agt}
{As already mentioned in the introduction, to unify the analysis of a large class of first-order methods, we will show how to construct an upper estimate $G^{(t)}$ of the optimality gap $f(\vxh^{(t)})-f(\vx^*)$, where $\vxh^{(t)}$ is the output of a first-order method at time $t$. This upper estimate is defined as $G^{(t)} = U^{(t)}-L^{(t)},$ where $U^{(t)}\geq f(\vxh^{(t)})$ is an upper bound on $f(\vxh^{(t)})$ and $L^{(t)}\leq f(\vx^*)$ is a lower bound on $f(\vx^*)$. We refer to $G^{(t)}$ as the approximate duality gap, due to the connections between the lower bound $L^{(t)}$ and the Fenchel dual of a certain approximation of the objective function $f(\vx^{(t)}),$ further discussed in Section~\ref{sec:lb}. To show that the method converges at some rate $\a^{(t)}$ (e.g., $\a^{(t)} = t$), we will show that $\a^{(t)}G^{(t)}$ is a non-increasing function of time, so that $\a^{(t)}G^{(t)}\leq \a^{(0)}G^{(0)},$ and, consequently, $f(\vxh^{(t)})-f(\vx^*)\leq G^{(t)}\leq \nicefrac{\a^{(0)}G^{(0)}}{\a^{(t)}}$. 
}

\subsection{Upper Bound} 
{The simplest upper bound on $f(\vxh^{(t)})$ is $f(\vxh^{(t)})$ itself: i.e., $U^{(t)} = f(\vxh^{(t)})$. In this case, $\vxh^{(t)}$ will be the last point constructed by the algorithm, i.e., $\vxh^{(t)} = \vx^{(t)}$. We will make this choice of the upper bound whenever we can assume that $f(\cdot)$ is differentiable (e.g., in the setting of accelerated and Frank-Wolfe methods), so that in the continuous time setting we can differentiate $\a^{(t)}U^{(t)}$ with respect to $t$ and write $\frac{\d}{\d t}(\a^{(t)}U^{(t)}) = \dot{\a}^{(t)}f(\vxh^{(t)}) + \a^{(t)}\innp{\nabla f(\vxh^{(t)}), \frac{\d}{\d t}{\vxh}^{(t)}}$. In the settings where $f(\cdot)$ is typically not assumed to be differentiable but only subdifferentiable (e.g., in the setting of dual averaging/mirror descent methods), $\vxh^{(t)}$ will be a weighted average of the points $\vx^{(\tau)}$ constructed by the method up to time $t$: $\vxh^{(t)} = \frac{\a^{(t)}-A^{(t)}}{\a^{(t)}}\vx^{(0)} + \frac{1}{\a^{(t)}}\int_{0}^t \vx^{(\tau)}\d\a^{(\tau)}$ and we will choose $U^{(t)} = \frac{\a^{(t)}-A^{(t)}}{\a^{(t)}}f(\vx^{(0)})+\frac{1}{\a^{(t)}}\int_0^t f(\vx^{(\tau)})\d\a^{(\tau)}$. Due to Jensen's inequality, $f(\vxh^{(t)})\leq U^{(t)},$ i.e., $U^{(t)}$ is a valid upper bound on $f(\vxh^{(t)})$. Observe that this choice of $U^{(t)}$ allows us to differentiate $\a^{(t)}U^{(t)}$ with respect to $t$ in the continuous time setting, and, thus, we can write $\frac{\d}{\d t}(\a^{(t)}U^{(t)}) = \dot{\a}^{(t)}f(\vx^{(t)}),$ as $\a^{(t)}-A^{(t)} = \a^{(0)}$ is a constant. These choices of upper bounds easily extend to the setting of composite 
objectives $\bar{f}(\cdot) = f(\cdot) + \psi(\cdot)$ (see Section~\ref{sec:ct-algos} for more details).}

\subsection{Lower Bound}\label{sec:lb}
{The simplest lower bound on $f(\vx^*)$ is $f(\vx^*)$ itself. However, it is not clear how to use $L^{(t)} = f(\vx^*)$ and guarantee $\frac{\d}{\d t}(\a^{(t)}G^{(t)})\leq 0,$ as in that case in the continuous-time domain $\frac{\d}{\d t}(\alpha^{(t)}L^{(t)}) = \dot{\a}^{(t)}f(\vx^*)$ is not possible to evaluate, as we do not know $f(\vx^*)$ (recall that, by assumption, $\dot{\a}^{(t)}>0$). Observe that if instead $f(\vx^*)$ appeared in the lower bound as $\frac{c}{\a^{(t)}}f(\vx^*)$ for some constant $c,$ we would not have this problem anymore, as $f(\vx^*)$ would not appear in $\frac{\d}{\d t}(\alpha^{(t)}L^{(t)})$. This is true because $\frac{\d}{\d t}\left(\a^{(t)}\frac{c}{\a^{(t)}}f(\vx^*)\right) = \frac{\d}{\d t}\left({c}f(\vx^*)\right)=0.$ On the other hand, convexity of $f(\cdot)$ leads to the following lower-bounding hyperplanes for all $\vx \in X$: $f(\vx^*) \geq f(\vx) + \innp{\nabla f(\vx), \vx^* - \vx}$.\footnote{{Observe here that if $f(\cdot)$ is not differentiable but only subdifferentiable, we can still obtain lower-bounding hyperplanes by using subgradients in place of the gradients.}} In particular, taking a convex combination of the trivial lower bound $f(\vx^*)$ and the lower-bounding hyperplanes defined by points $\vx^{(\tau)}$ constructed by the method up to time $t$, we have:}
\begin{equation}\label{eq:lb-simple}
{f(\vx^*) \geq \frac{\a^{(t)}-A^{(t)}}{\a^{(t)}}f(\vx^*) + \frac{1}{\a^{(t)}}\int_0^t \left(f(\vx^{(\tau)}) + \innp{\nabla f(\vx^{(\tau)}), \vx^*-\vx^{(\tau)}}\right)\d\a^{(\tau)}}.
\end{equation}
{
As in the continuous-time domain $\a^{(t)}-A^{(t)} = \a^{(0)}$ is a positive constant, the lower bound equal to the right-hand side of~\eqref{eq:lb-simple} is well-defined at $t=0$ and $f(\vx^*)$ does not appear in  $\frac{\d}{\d t}(\alpha^{(t)}L^{(t)})$. However, it would still not be possible to evaluate $\frac{\d}{\d t}(\alpha^{(t)}L^{(t)})$, as $\vx^*$ is not known. One way of addressing this issue is to use that:
\begin{equation}\label{eq:ineq-for-lb-fw}
\begin{aligned}
&\int_0^t \left(f(\vx^{(\tau)}) + \innp{\nabla f(\vx^{(\tau)}), \vx^*-\vx^{(\tau)}}\right)\d\a^{(\tau)}\\
&\hspace{2cm}\geq \int_0^t \min_{\vu \in X}\left(f(\vx^{(\tau)}) + \innp{\nabla f(\vx^{(\tau)}), \vu-\vx^{(\tau)}}\right)\d\a^{(\tau)}, \text{ or }
\end{aligned}
\end{equation}
\begin{equation}\label{eq:ineq-for-lb-gen}
\begin{aligned}
&\int_0^t \left(f(\vx^{(\tau)}) + \innp{\nabla f(\vx^{(\tau)}), \vx^*-\vx^{(\tau)}}\right)\d\a^{(\tau)}\\
&\hspace{2cm}\geq \min_{\vu \in X}\int_0^t \left(f(\vx^{(\tau)}) + \innp{\nabla f(\vx^{(\tau)}), \vu-\vx^{(\tau)}}\right)\d\a^{(\tau)}.
\end{aligned}
\end{equation}
While the inequality~\eqref{eq:ineq-for-lb-gen} is tighter than~\eqref{eq:ineq-for-lb-fw}, as a minimum of affine functions it is not differentiable w.r.t.~$\vu$ (and consequently not differentiable w.r.t.~$t$). The use of~\eqref{eq:ineq-for-lb-fw} leads to the continuous-time version of the standard Frank-Wolfe method~\cite{frank1956algorithm}.
\begin{example}\label{ex:frank-wolfe}
Standard continuous-time Frank-Wolfe method. Using~\eqref{eq:ineq-for-lb-fw}, we have the following lower bound:
\begin{equation}\label{eq:fw-lb}
\begin{aligned}
f(\vx^*) \geq L^{(t)} \defeq &\frac{\int_0^t \min_{\vu \in X}\left(f(\vx^{(\tau)}) + \innp{\nabla f(\vx^{(\tau)}), \vu-\vx^{(\tau)}}\right)\d\a^{(\tau)}}{\a^{(t)}}\\
&+ \frac{(\a^{(t)}-A^{(t)})f(\vx^*)}{\a^{(t)}}.
\end{aligned}
\end{equation}
Since the standard assumption in this setting is that $f(\cdot)$ is smooth (or, at the very least, continuously differentiable), we take $U^{(t)} = f(\vx^{(t)})$ and $\vxh^{(t)} = \vx^{(t)}$. Denoting $\vv^{(t)} \in \argmin_{\vu \in X}\left(f(\vx^{(\tau)}) + \innp{\nabla f(\vx^{(\tau)}), \vu-\vx^{(\tau)}}\right)$ and computing $\a^{(t)}G^{(t)},$ we get:
\begin{align*}
\frac{\d}{\d t}(\a^{(t)}G^{(t)}) =& \innp{\nabla f(\vx^{(t)}), \a^{(t)}\dot{\vx}^{(t)} - \dot{\a}^{(t)}(\vv^{(t)}-\vx^{(t)})}.
\end{align*}
Setting $\a^{(t)}\dot{\vx}^{(t)} - \dot{\a}^{(t)}(\vv^{(t)}-\vx^{(t)})=0$ gives $\frac{\d}{\d t}(\a^{(t)}G^{(t)}) = 0$ and precisely recovers the continuous-time version of the Frank-Wolfe algorithm, as in that case $\dot{\vx}^{(t)} = \frac{\dot{\a}^{(t)}(\vv^{(t)}-\vx^{(t)})}{\a^{(t)}},$ i.e., (as explained in Section~\ref{sec:prelims}) $\vx^{(t)}$ is a weighted average of $\vv^{(t)}$'s. 
\end{example}
Notice that the use of~\eqref{eq:ineq-for-lb-fw} in the construction of the lower bound makes sense only when linear minimization over $X$ is possible. However, there are insights we can take from the construction of the lower bound~\eqref{eq:fw-lb}. In particular, we can alternatively view~\eqref{eq:fw-lb} as being constructed as a lower bound on $f(\vx^*) + \psi(\vx^*),$ where $\psi(\vx^*)$ is the indicator of $X$. Hence we can view $L^{(t)}$ from~\eqref{eq:fw-lb} as being constructed as follows:
\begin{align*}
f(\vx^*) + \psi(\vx^*) \geq &\frac{\int_0^t \left(f(\vx^{(\tau)}) + \innp{\nabla f(\vx^{(\tau)}), \vx^*-\vx^{(\tau)}} + \psi(\vx^*)\right)\d\a^{(\tau)}}{\a^{(t)}}\\
&+ \frac{(\a^{(t)}-A^{(t)})f(\vx^*)}{\a^{(t)}}\geq L^{(t)}.
\end{align*}
We will see later in Section~\ref{sec:ct-algos} how this leads to a more general version of Frank-Wolfe method for composite functions, along the lines of the method from~\cite{nesterov2015cgm}.}

{Constructing a lower bound on $f(\vx^*) + \psi(\vx^*)$ when $\psi(\vx^*)$ was an indicator function had no effect, as in that case $f(\vx^*) + \psi(\vx^*) = f(\vx^*)$. To generalize this idea, we can construct a lower bound on a function that closely approximates $f(\cdot)$ around $\vx^*$. Since we want to obtain convergent methods, any error we introduce into this approximation should vanish at rate $1/\a^{(t)}$. Hence, a natural choice is to create a lower bound on $f(\vx^*) + \frac{1}{\a^{(t)}}\phi(\vx^*)$, where $\phi(\vx^*)$ is bounded\footnote{{A common choice of $\phi(\cdot)$ that ensures boundedness of $\phi(\vx^*)$ and non-negativity of $\phi(\cdot)$ is Bregman divergence of some function $\psi;$ namely, $\phi(\cdot) = D_{\psi}(\cdot, \vx^{(0)})$. Hence, in this case  we can interpret $\phi(\vx^*)$ as a generalized notion of the initial distance to the optimal solution.}}:
\begin{align*}
f(\vx^*) + \frac{1}{\a^{(t)}}\phi(\vx^*) \geq &\frac{\int_0^t f(\vx^{(\tau)})\d\a^{(\tau)} + \int_0^t \innp{\nabla f(\vx^{(\tau)}), \vx^* - \vx^{(\tau)}}\d\a^{(\tau)} + \phi(\vx^*)}{\a^{(t)}}\\
&+ \frac{(\a^{(t)}-A^{(t)})f(\vx^*)}{\a^{(t)}}.
\end{align*}
Now, if $\phi(\cdot)$ is strictly convex, $\min_{\vu \in X}\{\int_0^t \innp{\nabla f(\vx^{(\tau)}), \vu - \vx^{(\tau)}}\d\a^{(\tau)} + \phi(\vu)\}$ is differentiable (i.e., we can generalize the stronger inequality from~\eqref{eq:ineq-for-lb-gen}). We can view this as regularization of the minimum from~\eqref{eq:ineq-for-lb-gen}, leading to the following lower bound:
\begin{equation}\label{eq:lb-general}
\begin{aligned}
 f(\vx^*) +& \frac{1}{\a^{(t)}}\phi(\vx^*) \geq L^{(t)} + \frac{\phi(\vx^*)}{\a^{(t)}}\\
 =& \frac{\int_0^t f(\vx^{(\tau)})\d\a^{(\tau)} + \min_{\vu \in X}\{\int_0^t \innp{\nabla f(\vx^{(\tau)}), \vu - \vx^{(\tau)}}\d\a^{(\tau)} + \phi(\vu)\}}{\a^{(t)}}\\
 &+ \frac{(\a^{(t)}-A^{(t)})f(\vx^*)}{\a^{(t)}}.
\end{aligned}
\end{equation} 
Another advantage of the lower bound from~\eqref{eq:lb-general} over the previous one is that for many feasible sets $X$ there are natural choices of $\phi$ for which the minimization inside the lower bound from~\eqref{eq:lb-general} is easily solvable, often in a closed form (see, e.g.,~\cite{ben2001lectures}).}

\paragraph{Dual View of the Lower Bound}
{An alternative view of the lower bound from~\eqref{eq:lb-general} is through the concept of Fenchel Duality, which is defined for the sum of two convex functions (or the difference of a convex and a concave function). In particular, the Fenchel dual of $f(\vx)+\phi_t(\vx)$ is defined as $-f^*(-\vu)-\phi_t^*(\vu)$ (see, e.g., Chapter~15.2 in~\cite{bauschke2011convex}). Let $\phi_t(\cdot) = \frac{1}{A^{(t)}}\phi(\cdot)$ and $\vu^{(t)} = - \frac{\int_0^t \nabla f(\vx^{(\tau)})\d\a^{(\tau)}}{A^{(t)}}$. Observe that the minimization problem from~\eqref{eq:lb-general} defines a convex conjugate of $\phi_t.$ Thus, we can equivalently write~\eqref{eq:lb-general} as:
\begin{align*}
& f(\vx^*) + \frac{A^{(t)}}{\a^{(t)}}\phi_t(\vx^*) \geq  L^{(t)} + \frac{A^{(t)}}{\a^{(t)}}\phi_t(\vx^*)\\
&\hspace{1cm} =  \frac{\int_0^t (f(\vx^{(\tau)})-\innp{\nabla f(\vx^{(\tau)}), \vx^{(\tau)}})\d\a^{(\tau)}}{\a^{(t)}} - \frac{A^{(t)}}{\a^{(t)}}\phi_t^*(\vu^{(t)}) + \frac{\a^{(t)}-A^{(t)}}{\a^{(t)}}f(\vx^*)\\
&\hspace{1cm} = \frac{-\int_0^t f^*(\nabla f(\vx^{(\tau)}))\d\a^{(\tau)}}{\a^{(t)}} - \frac{A^{(t)}}{\a^{(t)}}\phi_t^*(\vu^{(t)}) + \frac{\a^{(t)}-A^{(t)}}{\a^{(t)}}f(\vx^*).
\end{align*}
Rearranging the terms in the last equality, we can equivalently write:
\begin{align*}
& L^{(t)} + \phi_t(\vx^*)\\
&\hspace{.5cm}= -\frac{A^{(t)}}{\a^{(t)}}\Big(\frac{\int_0^t f^*(\nabla f(\vx^{(\tau)}))\d\a^{(\tau)}}{A^{(t)}} + \phi_t^*(\vu^{(t)})\Big) + \frac{\a^{(t)}-A^{(t)}}{\a^{(t)}}(f(\vx^*) + \phi_t(\vx^*)).
\end{align*}
By Jensen's inequality, 
$
    \frac{\int_0^t f^*(\nabla f(\vx^{(\tau)}))\d\a^{(\tau)}}{A^{(t)}} \geq f^*(-\vu^{(t)}).
$\footnote{An earlier version of the paper made an incorrect use of Jensen's inequality at this point, and has since been corrected. We thank John Peebles for pointing it out.} 
Hence, we can view the general lower bound from~\eqref{eq:lb-general} as being slightly weaker than the weighted average of $f(\vx^*) + \phi_t(\vx^*)$ and its Fenchel dual $-f^*(-\vu^{(t)})-\phi_t^*(\vu^{(t)})$ evaluated at the average negative gradient $\vu^{(t)},$ and corrected by the introduced approximation error $\phi_t(\vx^*)$. Alternatively, the lower bound can be derived directly from $-f^*(-\vu^{(t)})-\phi_t^*(\vu^{(t)})$, slightly relaxing $-f^*(-\vu^{(t)})$ using the Jensen inequality. This means that we can roughly think about the lower bound $L^{(t)}$ as encoding the Fenchel dual of $f(\vx^*) + \phi_t(\vx^*)$ -- an approximation to $f(\vx^*)$ that converges to $f(\vx^*)$ at rate $1/\a^{(t)}$ -- and constructing dual solutions from the history of the gradients of $f$.}

\paragraph{Extension to Strongly Convex Objectives} When the objective is $\sigma$-strongly convex for some $\sigma > 0$, we can use $\sigma$-strong convexity (instead of just regular convexity) in the construction of the lower bound. This will generally give us a better lower bound which will  lead to the better convergence guarantees in the discrete-time domain. It is not hard to verify (by repeating the same arguments as above) that in this case we have the following lower bound:
\begin{equation}\label{eq:lb-strongly-cvx}
\begin{aligned}
L^{(t)} =&  \frac{\int_{0}^t f(\vx^{(\tau)}) \mathrm{d}\alpha^{(\tau)}}{\alpha^{(t)}}+ \frac{( \alpha^{(t)} - A^{(t)})f(\vx^*)  - \phi(\vx^*)}{\alpha^{(t)}}\\
&+ \frac{\min_{\vx \in X} \left\{\int_{0}^t \left(\innp{\nabla f(\vx^{(\tau)}), \vx - \vx^{(\tau)}} + \frac{\sigma}{2}\|\vx - \vx^{(\tau)}\|^2\right)\mathrm{d}\alpha { + \phi(\vx)}\right\}}{\alpha^{(t)}}.
\end{aligned}
\end{equation}
\begin{remark}
{Note that, due to the strong convexity of $f$, we do not need additional regularization in the lower bound to ensure that the minimum inside it is differentiable, i.e., we could use $\phi(\cdot) = 0$. This choice of $\phi$ will have no effect in the continuous-time convergence. In the discrete time, however, if we chose $\phi(\cdot) = 0$ the initial gap (and, consequently, the convergence bound) would scale with $\|\vx^{(1)} - \vx^{(0)}\|^2$. Adding a little bit of regularization (i.e., choosing a non-zero $\phi$) will allow us to replace $\|\vx^{(1)} - \vx^{(0)}\|^2$ with $\|\vx^{*} - \vx^{(0)}\|^2$.} 
\end{remark}
\paragraph{Extension to Composite Objectives} Suppose now that we have a composite objective $\bar{f}(\vx) = f(\vx) + \psi(\vx)$. Then, we can choose to apply the convexity argument only to $f(\cdot)$ and use $\psi(\cdot)$ as a regularizer (this will be particularly useful in the discrete-time domain in the settings where $f(\cdot)$ has some smoothness properties while $\psi(\cdot)$ is generally non-smooth). Therefore, we could start with $\bar{f}(\vx) \geq f(\vxh) + \innp{\nabla f(\vxh), \vx - \vxh} + \psi(\vx)$. Repeating the same arguments as in the general construction of the lower bound presented earlier in this subsection:
\begin{equation}\label{eq:lb-composite}
\begin{aligned}
L^{(t)} \defeq &\frac{\int_{0}^t f(\vx^{(\tau)}) \mathrm{d}\alpha^{(\tau)}}{\alpha^{(t)}}+ \frac{(\alpha^{(t)}-A^{(t)})\bar{f}(\vx^*) - \phi(\vx^*)}{\alpha^{(t)}}\\ 
&+ \frac{\min_{\vx \in X} \left\{\int_{0}^t \innp{\nabla f(\vx^{(\tau)}), \vx - \vx^{(\tau)}}\mathrm{d}\alpha^{(\tau)}  + A^{(t)}\psi(\vx) + \phi(\vx)\right\}}{\alpha^{(t)}}.
\end{aligned}
\end{equation}

\subsection{Extension to Monotone Operators and Saddle-Point Formulations} The notion of the approximate gap can be defined for problem classes beyond convex minimization. Examples are monotone operators and convex-concave saddle-point problems. {More details are provided in Appendix~\ref{app:mon-ops}}. 

\section{First-Order Methods in Continuous Time}\label{sec:ct-algos}

We now show how different {assumptions about the problem leading to the different} choices of the upper and lower bounds (and, consequently, the gap) yield different first-order methods. 

\subsection{Mirror Descent/Dual Averaging Methods}\label{sec:ct-md}

{Let us start by making minimal assumptions about the objective function $f$: we will assume that $f$ is convex and subdifferentiable (with the abuse of notation, in this case $\nabla f(\vx^{(t)})$ denotes an arbitrary but fixed subgradient of $f$ at $\vx^{(t)}$). As discussed in the previous section, since we are not assuming that $f$ is differentiable, we will take $U^{(t)} = \frac{\a^{(0)}}{\a^{(t)}}f(\vx^{(0)}) + \frac{\int_0^t f(\vx^{(\tau)})\d\a^{(\tau)}}{\a^{(t)}},$ so that $\a^{(t)}U^{(t)}$ is differentiable w.r.t.~the time and well-defined at the initial time point. As there are no additional assumptions about the objective (such as composite structure and strong convexity), we will use the generic lower bound from~\eqref{eq:lb-general}. Hence, we have the following expression for the gap:} 
\begin{equation}
\begin{aligned}
G^{(t)} = &\frac{-\min_{\vx \in X}\left\{ \int_{0}^t\innp{\nabla f(\vx^{(\tau)}), \vx - \vx^{(\tau)}}\mathrm{d}\alpha^{(\tau)} + \phi(\vx) \right\}}{\alpha^{(t)}}\\
&+ \frac{\alpha^{(0)}(f(\vx^{(0)})-f(\vx^*)) + \phi(\vx^*) }{\alpha^{(t)}}.
\end{aligned}
\end{equation}
Observe that $G^{(0)} \leq \frac{\phi(\vx^*)}{\alpha^{(0)}} + f(\vx^{(0)}) - f(\vx^*)$. Thus, if we show that $\frac{\d}{\d t}(\alpha^{(t)}G^{(t)})\leq 0$, this would immediately imply:
$$
f(\vxh^{(t)}) - f(\vx^*) \leq U^{(t)} - L^{(t)} \leq \frac{\alpha^{(0)}}{\alpha^{(t)}}G^{(0)} \leq \frac{\alpha^{(0)}(f(\vx^{(0)})-f(\vx^*)) + \phi(\vx^*)}{\alpha^{(t)}}.
$$

Now we show that ensuring the invariance $\frac{\d}{\d t}(\alpha^{(t)}G^{(t)}) = 0$ exactly produces the continuous-time mirror descent dynamics. Using (\ref{eq:mirror-map-z-t}) and Proposition (\ref{prop:ct-differentiation-inside-min}) with $\vz^{(t)} = -\int_{0}^t\nabla f(\vx^{(\tau)})\mathrm{d}\alpha^{(t)}$:
\begin{equation*}
\frac{\d}{\d t}(\alpha^{(t)}G^{(t)}) = -\innp{\nabla f(\vx^{(t)}), \nabla\phi^*(\vz^{(t)})-\vx^{(t)}}\dot{\alpha}^{(t)}.
\end{equation*}
Thus, to have $\frac{\d}{\d t}(\alpha^{(t)}G^{(t)}) = 0$, we can set $\vx^{(t)} = \nabla \phi^*(\vz^{(t)})$, which is precisely the mirror descent dynamics from \cite{nemirovskii1983problem}:
\begin{equation}\label{eq:ct-md}\tag{CT-MD}
\begin{gathered}
\dot{\vz}^{(t)} = - \dot{\alpha}^{(t)}\nabla f(\vx^{(t)}), \\
\vx^{(t)} = \nabla\phi^*(\vz^{(t)}), \\
\dot{\hat{\mathbf{x}}}^{(t)} = \dot{\alpha}^{(t)}\frac{\vx^{(t)}-\vxh^{(t)}}{\alpha^{(t)}},\\
\vz^{(0)} = 0,\; \vxh^{(0)} = \vx^{(0)}, \text{ for an arbitrary initial point } \vx^{(0)}\in X.
\end{gathered}
\end{equation}
{We immediately obtain the following lemma:}
\begin{lemma}\label{lemma:ct-md-conv}
Let $\vx^{(t)}, \vxh^{(t)}$ evolve according to (\ref{eq:ct-md}) for some convex function $f:X\rightarrow\mathbb{R}$. Then, $\forall t > 0$:
$$
f(\vxh^{(t)}) - f(\vx^*) \leq \frac{\alpha^{(0)}(f(\vx^{(0)})-f(\vx^*)) + \phi(\vx^*)}{\alpha^{(t)}}.
$$
\end{lemma}

\subsection{Accelerated Convex Minimization}\label{sec:ct-acc}
{Let us now assume more about the objective function: we will assume that $f$ is continuously differentiable, which means that our choice of the upper bound will be $U^{(t)} = f(\vxh^{(t)}) = f(\vx^{(t)}).$ Since there are no additional assumptions about $f,$ the lower bound we use is the generic one from Equation~\eqref{eq:lb-general}. Differentiating $\a^{(t)}G^{(t)},$ we have:}
\begin{align*}
\frac{\d}{\d t}(\alpha^{(t)}G^{(t)}) =& \frac{\d}{\d t}(\alpha^{(t)}f(\vx^{(t)})) - \dot{\alpha}^{(t)}\left(f(\vx^{(t)}) + \innp{\nabla f(\vx^{(t)}), \nabla \phi^*(\vz^{(t)})-\vx^{(t)}}\right)\\
=& \innp{\nabla f(\vx^{(t)}), \alpha^{(t)}\dot{\vx}^{(t)}- \dot{\alpha}^{(t)}(\nabla \phi^*(\vz^{(t)})-\vx^{(t)})},
\end{align*}
where we have used $\frac{\d}{\d t}(f(\vx^{(t)})) = \innp{\nabla f(\vx^{(t)}), \dot{\vx}^{(t)}}$. Choosing $\dot{\vx}^{(t)} = \dot{\alpha}^{(t)}\frac{\nabla \phi^*(\vz^{(t)})-\vx^{(t)}}{\alpha^{(t)}}$, we get $\frac{\d}{\d t}(\alpha^{(t)}G^{(t)}) = 0$. This is precisely the accelerated mirror descent dynamics \cite{krichene2015accelerated,wibisono2016variational}, and we immediately get the convergence result stated as Lemma \ref{lemma:ct-amd-conv} below.
\begin{equation}\label{eq:ct-amd}\tag{CT-AMD}
\begin{gathered}
\dot{\vz}^{(t)} = -\dot{\alpha}^{(t)}\nabla f(\vx^{(t)}),\\
\dot{\vx}^{(t)} = \dot{\alpha}^{(t)}\frac{\nabla \phi^*(\vz^{(t)})-\vx^{(t)}}{\alpha^{(t)}},\\
\vz^{(0)} = 0,\; \vx^{(0)}\in X \text{ is an arbitrary initial point}.
\end{gathered}
\end{equation}
\begin{lemma}\label{lemma:ct-amd-conv}
Let $\vx^{(t)}, \vz^{(t)}$ evolve according to (\ref{eq:ct-amd}), for some {continuously differentiable} convex function $f:X\rightarrow\mathbb{R}$. Then, $\forall t > 0$:
$$
f(\vx^{(t)})-f(\vx^*) \leq \frac{\alpha^{(0)}(f(\vx^{(0)})-f(\vx^*)) + \phi(\vx^*)}{\alpha^{(t)}}.
$$
\end{lemma}
\subsection{Gradient Descent}

Using the same approximate gap as in the previous subsection, now consider the special case when $\phi(\vx) = \frac{\sigma}{2}\|\vx - \vx^{(0)}\|_2^2$ for an arbitrary initial point $\vx^{(0)} = \vx^{(0)}\in X$, $X = \mathbb{R}^n$, and some $\sigma > 0$. Then, $\nabla \phi^*(\vz^{(t)}) = \vx^{(0)} + \nicefrac{\vz^{(t)}}{\sigma}$. Using the result for $\frac{\d}{\d t}(\alpha^{(t)}G^{(t)})$ from the previous subsection and setting $\vx^{(t)} = \nabla\phi^*(\vz^{(t)})$, we have:
\begin{align*}
\frac{\d}{\d t}(\alpha^{(t)}G^{(t)}) &= \alpha^{(t)}\innp{\nabla f(\vx^{(t)}), \dot{\vx}^{(t)}} = -\frac{\alpha^{(t)}\dot{\alpha}^{(t)}}{\sigma}\|\nabla f(\vx^{(t)})\|_2^2 \leq 0.
\end{align*}
The choice $\vx^{(t)} = \nabla\phi^*(\vz^{(t)}) = \vx^{(0)} + \nicefrac{\vz^{(t)}}{\sigma}$  precisely defines the gradient descent algorithm, and the convergence result stated as Lemma~\ref{lemma:ct-grad-desc} immediately follows.
\begin{equation}\label{eq:ct-gd}\tag{CT-GD}
\begin{gathered}
\dot{\vz}^{(t)} = - \dot{\alpha}^{(t)}\nabla f(\vx^{(t)})\\
\vx^{(t)}= \nabla \phi^*(\vz^{(t)}) = \vx^{(0)} + \nicefrac{\vz^{(t)}}{\sigma},\\
\vz^{(0)} = 0, \; \vx^{(0)}\in \mathbb{R}^n \text{ is an arbitrary initial point}.
\end{gathered}
\end{equation}
\begin{lemma}\label{lemma:ct-grad-desc}
Let $\vx^{(t)}, \vz^{(t)}$ evolve according to (\ref{eq:ct-gd}), for some {continuously differentiable} convex function $f:\mathbb{R}^n\rightarrow\mathbb{R}$. Then, $\forall t > 0$:
$$
f(\vx^{(t)})-f(\vx^*) \leq \frac{\alpha^{(0)}(f(\vx^{(0)})-f(\vx^*)) + \frac{\sigma}{2}\|\vx^* - \vx^{(0)}\|_2^2}{\alpha^{(t)}}.
$$
\end{lemma}
\subsection{Accelerated Strongly Convex Minimization}
{Let us now assume that, in addition to being continuously differentiable, $f$ is strongly convex. }
 In that case, we use $U^{(t)} = f(\vx^{(t)})$ and the lower bound from (\ref{eq:lb-strongly-cvx}). Let $\phi_t(\vx) = \int_{0}^t \frac{\sigma}{2}\|\vx - \vx^{(\tau)}\|^2\mathrm{d}\a^{(\tau)}{ + \phi(\vx)}$. Observe that $\frac{\d}{\d t}\phi_t(\vx) \geq 0$, $\forall \vx \in X$. Then, we have the following result for the change in the gap:
\begin{align*}
\frac{\d}{\d t}(\alpha^{(t)}G^{(t)}) =& \frac{\d}{\d t}(\alpha^{(t)}f(\vx^{(t)})) - \dot{\alpha}^{(t)}\left(f(\vx^{(t)}) + \innp{\nabla f(\vx^{(t)}), \nabla \phi_t^*(\vz^{(t)})-\vx^{(t)}}\right)\\
&- \frac{\d}{\d\tau}\left.(\phi_{\tau}(\nabla \phi_t^*(\vz^{(t)})))\right|_{\tau = t}\\
\leq& \innp{\nabla f(\vx^{(t)}), \alpha^{(t)}\dot{\vx}^{(t)}- \dot{\alpha}^{(t)}(\nabla \phi_t^*(\vz^{(t)})-\vx^{(t)})}.
\end{align*}
Therefore, choosing $\dot{\vx}^{(t)} = \dot{\alpha}^{(t)}\frac{\nabla \phi_t^*(\vz^{(t)})-\vx^{(t)}}{\alpha^{(t)}}$ gives $\frac{\d}{\d t}(\alpha^{(t)}G^{(t)})\leq 0$, and the convergence result stated as Lemma \ref{lemma:ct-amd-sc-conv} below follows.
\begin{equation}\label{eq:ct-amd-sc}\tag{CT-ASC}
\begin{gathered}
\dot{\vz}^{(t)} = -\dot{\alpha}^{(t)}\nabla f(\vx^{(t)}),\\
\dot{\vx}^{(t)} = \dot{\alpha}^{(t)}\frac{\nabla \phi_t^*(\vz^{(t)})-\vx^{(t)}}{\alpha^{(t)}},\\
\vz^{(0)} = 0,\; \vx^{(0)}\in \mathbb{R}^n \text{ is an arbitrary initial point}.
\end{gathered}
\end{equation}
\begin{lemma}\label{lemma:ct-amd-sc-conv}
Let $\vx^{(t)}$ evolve according to (\ref{eq:ct-amd-sc}) {for some continuously differentiable and $\sigma-$strongly convex function $F$}, where $\phi_t(\vx) = \int_{0}^t \frac{\sigma}{2}\|\vx - \vx^{(\tau)}\|^2{ + \phi(\vx)}$. Then, $\forall t > 0$:
$$
f(\vx^{(t)})-f(\vx^*) \leq \frac{\alpha^{(0)}(f(\vx^{(0)})-f(\vx^*)) + \phi(\vx^*)}{\alpha^{(t)}}.
$$
\end{lemma}
Note that, while there is no difference in the convergence bound for (\ref{eq:ct-amd}) and (\ref{eq:ct-amd-sc}), in the discrete time domain these two algorithms lead to very different convergence bounds, due to the different discretization errors they incur.
\subsection{Composite Dual Averaging}

Now assume that the objective is composite: $\bar{f}(\vx) = f(\vx) + \psi(\vx)$, where $f(\vx)$ is convex and $\nabla\phi_t^*(\cdot)$ is easily computable, for $\phi_t(\vx) = A^{(t)}\psi(\vx) + \phi(\vx)$. Then, we can use the lower bound for composite functions~(\ref{eq:lb-composite}). {Since we are not assuming that either $f$ or $\psi$ is continuously differentiable, the upper bound of choice is:} 
\begin{align*}
U^{(t)} &= \frac{1}{\alpha^{(t)}}\int_{0}^t \bar{f}(\vx^{(\tau)})\mathrm{d}\alpha^{(\tau)} + \frac{\alpha^{(0)}}{\alpha^{(t)}}\bar{f}(\vx^{(0)})\\
&= \frac{1}{\alpha^{(t)}}\int_{0}^t (f(\vx^{(\tau)})+\psi(\vx^{(\tau)}))\mathrm{d}\alpha^{(\tau)} + \frac{\alpha^{(0)}}{\alpha^{(t)}}(f(\vx^{(0)})+\psi(\vx^{(0)})).
\end{align*}
Then, the change in the gap is:
\begin{equation}\label{eq:ct-gap-change-comp}
\begin{aligned}
&\frac{\d}{\d t}(\alpha^{(t)}G^{(t)})\\
&\hspace{1cm}= \dot{\alpha}^{(t)}\psi(\vx^{(t)}) - \dot{\alpha}^{(t)}\innp{\nabla f(\vx^{(t)}), \nabla \phi_t^*(\vz^{(t)})-\vx^{(t)}} - \dot{\alpha}^{(t)}\psi(\nabla\phi_t^*(\vz^{(t)})).
\end{aligned}
\end{equation}
Thus, when $\dot{\vx}^{(t)} = \nabla\phi_t^*(\vz^{(t)})$, where $\phi_t(\cdot) = \alpha^{(t)}\psi(\cdot)$, we have $\frac{\d}{\d t}(\alpha^{(t)}G^{(t)}) = 0$, and Lemma \ref{lemma:ct-composite-conv} follows immediately. The algorithm can be thought of as mirror descent (or dual averaging) for composite minimization.
\begin{equation}\label{eq:ct-comp-md}\tag{CT-CMD}
\begin{gathered}
\dot{\vz}^{(t)} = -\dot{\alpha}^{(t)}\nabla f(\vx^{(t)}),\\
\vx^{(t)} = \nabla\phi_t^*(\vz^{(t)}),\\
\vxh^{(t)} = \dot{\alpha}^{(t)}\frac{\vx^{(t)}-\vxh^{(t)}}{\alpha^{(t)}},\\
\vz^{(0)} = 0,\; \vxh^{(0)} = \vx^{(0)}, \text{ for arbitrary initial point } \vx^{(0)}\in X.
\end{gathered}
\end{equation}
\begin{lemma}\label{lemma:ct-composite-conv}
Let $\vx^{(t)}, \vxh^{(t)}$ evolve according to (\ref{eq:ct-comp-md}), for some convex composite function $\bar{f} = f+ \psi: X\rightarrow \mathbb{R}$. Then, $\forall t > 0$:
$$
\bar{f}(\vxh^{(t)})-\bar{f}(\vx^*) \leq \frac{\alpha^{(0)}(\bar{f}(\vx^{(0)})-\bar{f}(\vx^*)) + \phi(\vx^*)}{\alpha^{(t)}}.
$$
\end{lemma}

\subsection{{Generalized} Frank-Wolfe Method} 
{We have already discussed the standard version of Frank-Wolfe method in Example~\ref{ex:frank-wolfe}. As discussed there, we can view standard Frank-Wolfe method as minimizing a composite objective $\bar{f} = f + \psi,$ where $\psi$ is the indicator function of set $X$, and using the assumption that  problems of the form $\min_{\vu \in X} \{\innp{\vz, \vu} + \psi(\vu)\}$ are easily solvable for any fixed $\vz$. We will now show how the method generalizes for any (possibly non-differentiable) $\psi$. The lower bound from~\eqref{eq:fw-lb} derived for standard Frank-Wolfe method immediately generalizes to:}
\begin{equation}\label{eq:gen-lb-fw}
\begin{aligned}
{
L^{(t)}} = &{\frac{\int_0^t f(\vx^{(\tau)})\d\a^{(\tau)} + \int_0^t \min_{\vu \in X}\{\innp{\nabla f(\vx^{(\tau)})\d\a^{(\tau)}, \vu - \vx^{(\tau)}} + \psi(\vu)\}}{\a^{(t)}}}\\
&{+ \frac{(\a^{(t)}-A^{(t)})\bar{f}(\vx^*)}{\a^{(t)}}.
}
\end{aligned}
\end{equation}
{By (generalized) Danskin's Theorem~\cite{bertsekas1971control}, we have that $$\nabla\psi^*(-\nabla f(\vx^{(\tau)}))\in \argmin_{\vu \in X}\Big\{\innp{\nabla f(\vx^{(\tau)})\d\a^{(\tau)}, \vu - \vx^{(\tau)}} + \psi(\vu)\Big\}$$ 
(note that the minimizer may not be unique as $\psi$ is not necessarily strictly convex; with the abuse of notation, $\nabla\psi^*(-\nabla f(\vx^{(\tau)}))$ may be a subgradient of $\psi^*$). Since $f$ is assumed to be continuously differentiable, we would like to use $f(\vx^{(t)})$ in the upper bound. On the other hand, $\psi$ is not necessarily differentiable, and hence we cannot use $\psi(\vx^{(t)})$ in the upper bound. Instead, we need to have $\frac{\a^{(0)}}{\a^{(t)}}\psi(\vx^{(0)}) + \frac{\int_0^t \psi(\vv^{(\tau)})\d\a^{(\tau)}}{\a^{(t)}}$ for some points $\vv^{(\tau)} \in X$ to ensure that $\a^{(t)}U^{(t)}$ is differentiable. Hence, based on the rules for choosing the upper bound discussed in Section~\ref{sec:agt}, we would like the upper bound to be of the form:}
\begin{equation}\label{eq:ub-fw}
{U^{(t)} = f(\vx^{(t)}) + \frac{\a^{(0)}}{\a^{(t)}}\psi(\vx^{(0)}) + \frac{\int_0^t \psi(\vv^{(\tau)})\d\a^{(\tau)}}{\a^{(t)}}.}
\end{equation}
{Using Jensen's Inequality, this is a valid upper bound on $\bar{f}(\vx^{(t)})$ if $\vx^{(t)} = \frac{\a^{(0)}}{\a^{(t)}}\vx^{(0)} + \frac{\int_0^t \vv^{(\tau)}\d\a^{(\tau)}}{\a^{(t)}}.$ Taking a leap of faith, let us consider the gap constructed based on the lower and upper bounds from~\eqref{eq:gen-lb-fw}, \eqref{eq:ub-fw}. Differentiating $\a^{(t)}G^{(t)}:$} 
\begin{align*}
{\frac{\d}{\d t}(\a^{(t)}G^{(t)}) =}&{\innp{\nabla f(\vx^{(t)}), \a^{(t)}\dot{\vx}^{(t)} - \dot{\a}^{(t)}(\nabla \psi^*(-\nabla f(\vx^{(t)}))-\vx^{(t)})}} \\
&{+ \dot{\a}^{(t)}\psi(\vv^{(t)}) - \dot{\a}^{(t)}\psi(\nabla \psi^*(-\nabla f(\vx^{(t)}))).}
\end{align*}
{For the terms from the second line to cancel out, we need $\vv^{(t)} = \nabla \psi^*(-\nabla f(\vx^{(t)})).$ Then, if we set $\a^{(t)}\dot{\vx}^{(t)} - \dot{\a}^{(t)}(\vv^{(t)}-\vx^{(t)}) = 0$, we get $\frac{\d}{\d t}(\a^{(t)}G^{(t)}) = 0$. But this precisely defines $\vx^{(t)}$ as $\vx^{(t)} = \frac{\a^{(0)}}{\a^{(t)}}\vx^{(0)} + \frac{\int_0^t \vv^{(\tau)}\d\a^{(\tau)}}{\a^{(t)}}$, which is what we needed for the upper bound to be valid. Hence, we exactly recover the continuous-time counterpart of the generalized Frank-Wolfe method from~\cite{nesterov2015cgm} and Lemma~\ref{lemma:ct-fw-conv} follows.}
\begin{equation}\label{eq:ct-fw}\tag{CT-FW}
\begin{gathered}
\vzh^{(t)} = -\nabla f(\vx^{(t)}),\\
\dot{\vx}^{(t)} = \dot{\alpha}^{(t)}\frac{\nabla\psi^*(\vzh^{(t)})-\vx^{(t)}}{\alpha^{(t)}},\\
\vx^{(0)}\in X \text{ is an arbitrary initial point}.
\end{gathered}
\end{equation}
\begin{lemma}\label{lemma:ct-fw-conv}
Let $\vx^{(t)}, \vzh^{(t)}$ evolve according to (\ref{eq:ct-fw}),  for some convex composite function $\bar{f} = f+ \psi: X\rightarrow \mathbb{R}${, where $f$ is continuously differentiable}. Then:
$$
\bar{f}(\vx^{(t)})-\bar{f}(\vx^*) \leq \frac{\alpha^{(0)}(\bar{f}(\vx^{(0)})-\bar{f}(\vx^*))}{\alpha^{(t)}}, \; \forall t \geq 0.
$$
\end{lemma}

\section{Discretization and Incurred Errors}\label{sec:discretization}

Suppose now that $\alpha^{(t)}$ is a discrete measure. In particular, let $\alpha^{(t)}$ be an increasing piecewise constant function, with $\alpha^{(t)}= 0$ for $t < 0$, $\alpha^{(t)}$ constant in intervals $(0 + i,\; 0 + i + 1)$ for $i \in \mathbb{Z}_+$, and $\alpha^{((0 + i)+)} - \alpha^{((0 + i)-)} = a_i$ for some $a_i > 0$ and $i \in \mathbb{Z}_+${, as discussed in Section~\ref{sec:notation}}.

For the continuous-time algorithms (and their analyses) presented in Section \ref{sec:ct-algos}, there are generally two causes of the discretization error: (i) different integration rules applying to continuous and discrete measures, and (ii) discontinuities in the algorithm updates. We discuss these two causes in more details below.  

\paragraph{Integration errors}
To understand where the integration errors occur, we first note that such errors cannot occur in integrals whose sole purpose is weighted averaging, since for these integrals there is no functional difference in the continuous- and discrete-time domains. 
Thus, the only place where the integration errors can occur is in the integral appearing under the minimum in the lower bound. In $\alpha^{(t)}G^{(t)} = A^{(t)}G^{(t)}$, the integral appears as:
\begin{equation}
I^{(0, t)} = -\int_{0}^t \innp{\nabla f(\vx^{(\tau)}), \nabla \phi_t^*(\vz^{(t)})-\vx^{(\tau)}}\mathrm{d}\alpha^{(\tau)},\notag
\end{equation}
where $\phi_t(\cdot) = \phi(\cdot)$ in the case of mirror descent and accelerated convex minimization. Let $I_c^{(0, t)}$ denote the value $I^{(0, t)}$  {would take if} $\alpha$ was a continuous measure, {i.e., if rules of continuous integration applied}. 

Observe that, as between times $i-1$ and $i$ $\dot{\alpha}^{(\tau)}$ samples the function under the integral at time $i$, we have:
\begin{equation}\label{eq:I-i-1-i}
I^{(i-1, i)} = -{a_i}\innp{\nabla f(\vx^{(i)}), \nabla\phi_i^*(\vz^{(i)}) - \vx^{(i)}}.
\end{equation}

\paragraph{Discontinuities in the algorithm updates} 

In all of the described algorithms, the updates for $\vx^{(t)}$ (and possibly $\vxh^{(t)}$) depend on $\nabla \phi_t^*(\vz^{(t)})$. Recall that $\vz^{(t)}$ aggregates negative gradients up to time $t$ and thus also depends on $\nabla f(\vx^{(t)})$. In the continuous time domain, this is not a problem, since the updates in $\vx^{(t)}$ can follow updates in $\vz^{(t)}$ with an arbitrarily small delay, meaning that in the limit we can take that $\vx^{(t)}$ changes simultaneously with $\vz^{(t)}$. In the discrete time, however, the delay between the two updates cannot be neglected, and {using implicit updates of the form} $\vx^{(t)} = g(\nabla \phi_t^*(\vz^{(t)}))$ for some function $g(\cdot)$ is in general either not possible or requires many fixed-point iterations. 

Apart from affecting the value of $I^{(i-1, i)}$ described above, the discontinuities will also contribute additional discretization error in the case of composite minimization. The reason for the additional discretization error is that the analysis of the gap reduction relies on bounding the change in $\psi(\vx^{(t)})$ (or the $\dot{\alpha}^{(\tau)}$-weighted average of $\psi(\vx^{(\tau)})$'s for $\tau \in [0, t]$) from the upper bound by $\dot{\alpha}^{(t)}\psi(\nabla\phi_t^*(\vz^{(t)}))$ from the lower bound. For composite dual averaging, this discretization error at time $i$ will amount to $a_i\left(\psi(\vx^{(i)})-\psi(\nabla\phi_i^*(\vz^{(i)}))\right)$. 
Similar to composite dual averaging, Frank-Wolfe will accrue discretization error   $a_i\left(\psi(\vx^{(i)})-\psi(\nabla\psi^*(\vzh^{(i)}))\right)$.

\paragraph{Effect of discretization errors on the gap}

Since in continuous time we had that $\frac{\d}{\d t}(\alpha^{(t)}G^{(t)})\leq 0$, if the discretization error between discrete time points $i-1$ and $i$ is $E_d^{(i)}$, then $A^{(i)}G^{(i)}-A^{(i-1)}G^{(i-1)}\leq E_d^{(i)}$, and we can conclude that:
\begin{equation}\label{eq:gap-w-discr-err}
G^{(k)} \leq \frac{a_0 G^{(0)}}{A^{(k)}} + \frac{\sum_{i=1}^k E_d^{(i)}}{A^{(k)}}.
\end{equation}
{Note that $E_d^{(i)}$'s contain both the integration error and the error due to the discontinuities in the algorithm updates discussed above. These discretization errors will ultimately determine the algorithms' convergence rate: while in the continuous time domain we could choose $\a^{(t)}$ (and $A^{(t)}$) to grow arbitrarily fast as a function of time, in the discrete time domain,  the discretization errors $E_d^{(k)}$ will depend on the choice of $A^{(k)}$. In particular, this co-dependence between $A^{(k)}$ and $E_d^{(k)}$ will determine the choice of $A^{(k)}$ leading to the highest decrease in the bound on $G^{(k)}$ from~\eqref{eq:gap-w-discr-err} as a function of $k$.}

We are now ready to bound the discretization errors of the algorithms from Section~\ref{sec:ct-algos}. Before doing so, we make the following two remarks: 
\begin{remark}The versions of mirror descent and mirror prox presented here are in fact the ``lazy'' versions of these methods that are known as Nesterov's dual averaging \cite{nesterov2009primal}.{The ``lazy'' and standard versions of the methods are equivalent whenever $X = \mathbb{R}^n$. The reason we chose to work with the ``lazy'' versions of the methods is that they follow more directly from the discretization of continuous-time dynamics presented in the previous section.} 
\end{remark}
\begin{remark}
It is possible to obtain the discretization error $E_d^{(i)}$ by directly computing $A^{(k)}G^{(k)}-A^{(k-1)}G^{(k-1)}$ for the discrete version of the gap. We have chosen the approach presented here to illustrate the effects of discretization.
\end{remark}

\subsection{Dual Averaging (Lazy Mirror Descent)}

Recall that in mirror descent, $\phi_i(\cdot) = \phi(\cdot)$. The discretization error can be bounded as follows. 

\begin{proposition}\label{prop:md-discr-err}
The discretization error for (\ref{eq:ct-md}) is:
\begin{equation}
E_d^{(i)} = -a_i\innp{\nabla f(\vx^{(i)}), \nabla\phi^*(\vz^{(i)})-\vx^{(i)}} - D_{\phi^*}(\vz^{(i-1)}, \vz^{(i)}).\notag
\end{equation}
\end{proposition}
\begin{proof}
In 
(\ref{eq:ct-md}), $\vx^{(\tau)}=\nabla\phi^*(\vz^{(\tau)})$, and thus:
\begin{align}
I_c^{(i-1, i)} &= \int_{i-1}^{i}\innp{\dot{\vz}^{(\tau)}, \nabla \phi^*(\vz^{(i)})-\nabla\phi^*(\vz^{(\tau)})}d\tau\notag\\
&= - \int_{i-1}^{i} \frac{d D_{\phi^*}(\vz^{(\tau)}, \vz^{(i)})}{d\tau}d\tau = D_{\phi^*}(\vz^{(i-1)}, \vz^{(i)}). \label{eq:md-I_c}
\end{align}
Since for non-composite functions $E_d^{(i)} = I^{(i-1, i)}-I_c^{(i-1, i)}$, combining (\ref{eq:I-i-1-i}) and (\ref{eq:md-I_c}) completes the proof.
\end{proof}

We now consider two different discretization methods that lead to discrete-time algorithms known as (lazy) mirror descent and mirror prox.

\paragraph{Forward Euler Discretization: Lazy Mirror Descent} 

Forward Euler discretization leads to the following algorithm updates:
\begin{equation}\label{eq:classical-md}\tag{MD}
\begin{gathered}
\vz^{(i)} = \vz^{(i-1)} - a_{i}\nabla f(\vx^{(i)}),\\
\vx^{(i)} = \nabla\phi^*(\vz^{(i-1)}),\\
\vxh^{(i)} = \frac{A^{(i-1)}}{A^{(i)}}\vxh^{(i-1)} + \frac{a_i}{A^{(i)}}\vx^{(i)},\\
\vz^{(0)} = -a_0 f(\vx^{(0)}), \; \vxh^{(0)} = \vx^{(0)}, \text{ for arbitrary } \vx^{(0)}\in X.
\end{gathered}
\end{equation}
It follows (from Proposition \ref{prop:md-discr-err}) that in this case the discretization error is given as:
\begin{equation}
E_d^{(i)} = -a_i\innp{\nabla f(\vx^{(i)}), \vx^{(i+1)}-\vx^{(i)}} - D_{\phi^*}(\vz^{(i-1)}, \vz^{(i)}).
\end{equation}
When $f(\cdot)$ is Lipschitz-continuous, we recover the classical mirror descent/dual averaging convergence result \cite{nemirovskii1983problem}:
\begin{theorem} \label{thm:md}
Let $f:X\rightarrow \mathbb{R}$ be an $L$-Lipschitz-continuous convex function, and let $\psi:X\rightarrow \mathbb{R}$ be $\sigma$-strongly convex for some $\sigma >0$. Let $\vx^{(i)}, \vxh^{(i)}$ evolve according to (\ref{eq:classical-md}) for $i\leq k$ and $k\geq 1$. Then, if $a_i = \frac{1}{L}\sqrt{\frac{2\sigma D_{\psi}(\vx^*, \vx^{(0)})}{k+1}}$ and 
$\phi(\cdot) = D_{\psi}(\cdot, \vx^{(0)})$:
$$
f(\vxh^{(k)}) - f(\vx^*) \leq \sqrt{\frac{2D_{\psi}(\vx^*, \vx^{(0)})}{\sigma}}\cdot \frac{L}{\sqrt{k+1}}.
$$
\end{theorem}
\begin{proof}
By Proposition \ref{prop:cvx-conj-bd-is-strongly-cvx-too}, $D_{\psi^*}(\vz^{(i-1)}, \vz^{(i)})\geq \frac{\sigma}{2}\|\nabla\psi^*(\vz^{(i-1)})-\nabla\psi^*(\vz^{(i)})\|^2 = \frac{\sigma}{2}\|\vx^{(i+1)}-\vx^{(i)}\|^2$. As $D_{\phi}(\cdot, \cdot) = D_{\psi}(\cdot, \cdot)$ and $f(\cdot)$ is Lipschitz continuous with parameter $L$, using Cauchy-Schwartz Inequality:
\begin{align*}
E_d^{(i)} \leq  a_i L\|\vx^{(i+1)}-\vx^{(i)}\| - \frac{\sigma}{2}\|\vx^{(i+1)}-\vx^{(i)}\|^2\leq \frac{{a_i}^2 L^2}{2\sigma},
\end{align*}
where the second inequality follows from $2ab - b^2 \leq a^2$, $\forall a, b$. Therefore, from (\ref{eq:gap-w-discr-err}):
\begin{equation}\label{eq:md-gap-w-error}
G^{(k)} \leq \frac{a_0 G^{(0)}}{A^{(k)}} + \frac{L^2}{2\sigma}\cdot \frac{\sum_{i=1}^k {a_i}^2}{A^{(k)}}.
\end{equation}
Similarly, we can bound the initial gap as:
\begin{align}
a_0 G^{(0)} &= -a_0\innp{\nabla f(\vx^{(0)}), \nabla \phi^*(\vz^{(0)})-\vx^{(0)}} - \phi(\nabla \phi^*(\vz^{(0)})) + \phi(\vx^*)\notag\\
&=  -a_0\innp{\nabla f(\vx^{(0)}), \vx^{(1)}-\vx^{(0)}} -D_{\psi}(\vx^{(1)}, \vx^{(0)}) + D_{\psi}(\vx^*, \vx^{(0)})\notag\\
&\leq \frac{{a_0}^2 L^2}{2\sigma} + D_{\psi}(\vx^*, \vx^{(0)}). \label{eq:md-initial-gap}
\end{align}
Finally, combining (\ref{eq:md-gap-w-error}), (\ref{eq:md-initial-gap}), the choice of $a_i$'s, and $f(\vxh^{(k)})-f(\vx^*)\leq G^{(k)}$, the result follows.
\end{proof}
\paragraph{Approximate Backward Euler Discretization: Mirror Prox/Extra-gradient} 

Observe that if we could set $\vx^{(i)}=\nabla\phi^*(\vz^{(i)})$ (i.e., if we were using backward Euler discretization for $\vx$), then the discretization error would be negative: $E_d^{(i)} = - D_{\phi^*}(\vz^{(i-1)}, \vz^{(i)})$. However, backward Euler is only an implicit discretization method, as it involves solving $\vx^{(i)} = \nabla\phi^*(\vz^{(i-1)}-a_i \nabla f(\vx^{(i)}))$. Fortunately, the fact that the discretization error is negative enables an approximate implementation of the method, where only two fixed-point iteration steps are performed.\footnote{{This discretization can also be viewed as the predictor-corrector method.}} The resulting discrete-time method is known as mirror prox \cite{Mirror-Prox-Nemirovski} or extra-gradient descent \cite{extragradient-descent}.
\begin{equation}\label{eq:mp}\tag{MP}
\begin{gathered}
\vxt^{(i-1)} = \nabla\phi^*(\vz^{(i-1)}),\\
\vzt^{(i-1)} = \vz^{(i-1)} - a_i \nabla f(\vxt^{(i-1)}),\\
\vx^{(i)} = \nabla\phi^*(\vzt^{(i-1)}),\\
\vz^{(i)} = \vz^{(i-1)} - a_i \nabla f(\vx^{(i)}), \\
\vxh^{(i)} = \frac{A^{(i-1)}}{A^{(i)}}\vxh^{(i)} + \frac{a_i}{A^{(i)}}\vx^{(i)},\\
\vz^{(0)} = - a_0 \nabla f(\vx^{(0)}), \text{ and } \vx^{(0)}\in X \text{ is an arbitrary initial point.}
\end{gathered}
\end{equation}
{This method is typically used for solving variational inequalities with monotone operators~\cite{Mirror-Prox-Nemirovski}. Its convergence bound is provided in Theorem~\ref{thm:mp-conv} in Appendix~\ref{app:mon-ops}.} 

\subsection{Accelerated Smooth Minimization}\label{sec:amd}

In this and in the following subsection, we will only consider  forward Euler discretization of the accelerated dynamics, which corresponds to the Nesterov's accelerated algorithm. Approximate backward Euler discretization using similar ideas as in the proof of convergence of mirror prox from the previous subsection is also possible and leads to the recent accelerated extra-gradient descent (\axgd) algorithm that we presented in \cite{AXGD}. 

As before, we can bound the discretization error by computing $I^{(i-1, i)}$ to obtain the following result.
\begin{proposition}\label{prop:amd-discr-err}
The discretization for (\ref{eq:ct-amd}) is:
\begin{equation}\label{eq:amd-fe-discr-error}
\begin{aligned}
E_d^{(i)} =& -a_i\innp{\nabla f(\vx^{(i)}), \nabla \phi^*(\vz^{(i)})-\vx^{(i)}} + A^{(i-1)}(f(\vx^{(i)})-f(\vx^{(i-1)}))\\
&- D_{\phi^*}(\vz^{(i-1)}, \vz^{(i)})\\
\leq & \innp{\nabla f(\vx^{(i)}), A^{(i)} \vx^{(i)} - A^{(i-1)}\vx^{(i-1)}-a_i \nabla \phi^*(\vz^{(i)})} - D_{\phi^*}(\vz^{(i-1)}, \vz^{(i)}).
\end{aligned}
\end{equation}
\end{proposition}
\begin{proof}
Recall continuous-time accelerated dynamics (\ref{eq:ct-amd}), where $\dot{\vx}^{(t)} = \dot{\alpha}^{(t)}\frac{\nabla\phi^*(\vz^{(t)})-\vx^{(t)}}{\alpha^{(t)}}$ and $\phi_i(\cdot) = \phi(\cdot)$. We have:
\begin{align*}
I_c^{(i-1, i)} &= -\int_{i-1}^{i}\innp{\nabla f(\vx^{(\tau)}), \nabla \phi^*(\vz^{(i)})-\vx^{(\tau)}}\mathrm{d}\alpha^{(\tau)}\\
&= -\int_{i-1}^{i}\innp{\nabla f(\vx^{(\tau)}), \alpha^{(\tau)}\dot{\vx}^{(\tau)}}d\tau + \int_{i-1}^{i}\innp{\dot{\vz}^{(\tau)}, \nabla \phi^*(\vz^{(i)}) - \nabla \phi^*(\vz^{(\tau)})}d\tau.
\end{align*}
Integrating by parts, the first integral is $-A^{(i-1)}(f(\vx^{(i)})-f(\vx^{(i-1)}))$, while the second integral is (as we have seen in the previous subsection) $D_{\psi^*}(\vz^{(i-1)}, \vz^{(i)})$. Thus, using (\ref{eq:I-i-1-i}), the discretization error is:
\begin{equation*}
\begin{aligned}
E_d^{(i)} =& -a_i\innp{\nabla f(\vx^{(i)}), \nabla \phi^*(\vz^{(i)})-\vx^{(i)}} + A^{(i-1)}(f(\vx^{(i)})-f(\vx^{(i-1)}))\\
&- D_{\phi^*}(\vz^{(i-1)}, \vz^{(i)})\\
\leq &\innp{\nabla f(\vx^{(i)}), A^{(i)} \vx^{(i)} - A^{(i-1)}\vx^{(i-1)}-a_i \nabla \phi^*(\vz^{(i)})} - D_{\phi^*}(\vz^{(i-1)}, \vz^{(i)}),
\end{aligned}
\end{equation*}
where we used $f(\vx^{(i)})-f(\vx^{(i-1)})\leq \innp{\nabla f(\vx^{(i)}), \vx^{(i)}-\vx^{(i-1)}}$, by $f(\cdot)$'s convexity. 
\end{proof}
Standard forward Euler discretization sets $\vx^{(i)} = \frac{A^{(i-1)}}{A^{(i)}}\vx^{(i-1)} + \frac{a_i}{A^{(i)}}\nabla \phi^*(\vz^{(i-1)})$, which results in the discretization error equal to $D_{\phi^*}(\vz^{(i)}, \vz^{(i-1)})$. We cannot bound such a discretization error, since we are not assuming that $f(\cdot)$ is Lipschitz-continuous. However, since $f(\cdot)$ is $L$-smooth, we can introduce an additional gradient step whose role is to cancel out the discretization error by reducing the upper bound. The algorithm then becomes the familiar Nesterov's accelerated method~\cite{Nesterov1983}:
\begin{equation}\label{eq:amd}\tag{AMD}
\begin{gathered}
\vz^{(i)} = \vz^{(i-1)} - a_i \nabla f(\vx^{(i)}),\\
\vx^{(i)} = \frac{A^{(i-1)}}{A^{(i)}}\vxh^{(i-1)} + \frac{a_i}{A^{(i)}}\nabla\phi^*(\vz^{(i-1)}),\\
\vxh^{(i)} = \mathrm{Grad}(\vx^{(i)}), \\
\vz^{(0)} = -a_0 f(\vx^{(0)}), \; \vxh^{(0)} = \mathrm{Grad}(\vx^{(0)}), \text{ for arbitrary } \vx^{(0)}\in X,
\end{gathered}
\end{equation}
where
\begin{equation}\label{eq:grad-step-smooth}
\mathrm{Grad}(\vx^{(i)}) = \arg\min_{\vx \in X}\left\{\innp{\nabla f(\vx^{(i)}), \vx - \vx^{(i)}} + \frac{L}{2}\|\vx - \vx^{(i)}\|^2\right\}. 
\end{equation}
The introduced gradient steps only affect the upper bound, changing it from $U^{(i)} = f(\vx^{(i)})$ to $U^{(i)} = f(\vxh^{(i)})$. Thus, correcting (\ref{eq:amd-fe-discr-error}) for the change in the upper bound:
\begin{align}
E_d^{(i)} =& -a_i\innp{\nabla f(\vx^{(i)}), \nabla \phi^*(\vz^{(i)})-\vx^{(i)}} + A^{(i-1)}(f(\vx^{(i)})-f(\vx^{(i-1)}))\notag\\
&- D_{\phi^*}(\vz^{(i-1)}, \vz^{(i)})
+ A^{(i)}(f(\vxh^{(i)}) - f(\vx^{(i)})) - A^{(i-1)}(f(\vxh^{(i-1)})-f(\vx^{(i-1)}))\notag\\
\leq & A^{(i)}(f(\vxh^{(i)}) - f(\vx^{(i)})) + \innp{\nabla f(\vx^{(i)}), A^{(i)} \vx^{(i)} - A^{(i-1)}\vxh^{(i-1)}-a_i \nabla \phi^*(\vz^{(i)})}\notag\\
&- D_{\phi^*}(\vz^{(i-1)}, \vz^{(i)})\notag\\
=& A^{(i)}(f(\vxh^{(i)}) - f(\vx^{(i)})) + a_i \innp{\nabla f(\vx^{(i)}), \nabla \phi^*(\vz^{(i-1)}) - \nabla \phi^*(\vz^{(i)})}\notag\\
&- D_{\phi^*}(\vz^{(i-1)}, \vz^{(i)}). \label{eq:amd-discr-err}
\end{align}

We are now ready to prove the convergence of Nesterov's algorithm for smooth functions \cite{Nesterov1983}:
\begin{theorem}\label{thm:amd-conv}
Let $f:X\rightarrow \mathbb{R}$ be an $L$-smooth function, $\psi:X\rightarrow \mathbb{R}$ be a $\sigma$-strongly convex function, and let $\phi(\cdot) = D_{\psi}(\cdot, \vx^{(0)})$. If $\vx^{(t)}, \vxh^{(t)}$ evolve according to (\ref{eq:amd}) for $a_i = {\frac{\sigma}{L}}\frac{i+1}{2}$, then $\forall k \geq 1$:
$$
f(\vxh^{(k)}) - f(\vx^*) \leq \frac{4L}{\sigma}\cdot\frac{D_{\psi}(\vx^*, \vx^{(0)})}{(k+1)(k+2)}.
$$
\end{theorem}
\begin{proof}
As $f(\cdot)$ is $L$-smooth, by the definition of $\vxh^{(i)}$:
\begin{equation}\label{eq:amd-f-smoothness}
f(\vxh^{(i)}) \leq f(\vx^{(i)}) + \min_{\vx \in X}\left\{\innp{\nabla f(\vx^{(i)}), \vx - \vx^{(i)}} + \frac{L}{2}\|\vx - \vx^{(i)}\|^2\right\}. 
\end{equation}
Since the gradient step was introduced to cancel out the discretization error, intuitively, it is natural to try to cancel out the second two terms from (\ref{eq:amd-discr-err}) (that correspond to the original discretization error (\ref{eq:amd-fe-discr-error})) by $A^{(i)}(f(\vxh^{(i)})-f(\vx^{(i)}))$, the decrease due to the gradient step. A point $\vx\in X$ that would charge the gradient term from (\ref{eq:amd-discr-err}) to the gradient term in (\ref{eq:amd-f-smoothness}) is $\vx = \vx^{(i)} - \frac{a_i}{A^{(i)}}(\nabla \phi^*(\vz^{(i-1)}) - \nabla \phi^*(\vz^{(i)})) = \frac{A^{(i-1)}}{A^{(i)}}\vxh^{(i-1)} + \frac{a_i}{A^{(i)}}\nabla \phi^*(\vz^{(i)}) \in X$. It follows from (\ref{eq:amd-f-smoothness}) that:
\begin{equation}\label{eq:amd-f-particular-x}
\begin{aligned}
A^{(i)}f(\vxh^{(i)}) \leq & A^{(i)}f(\vx^{(i)}) - {a_i}\innp{\nabla f(\vx^{(i)}), \nabla \phi^*(\vz^{(i-1)}) - \nabla \phi^*(\vz^{(i)})}\\
&+ \frac{L {a_i}^2}{2A^{(i)}}\|\nabla \phi^*(\vz^{(i-1)}) - \nabla \phi^*(\vz^{(i)})\|^2.
\end{aligned}
\end{equation}

By Proposition \ref{prop:cvx-conj-bd-is-strongly-cvx-too}, $D_{\phi^*}(\vz^{(i-1)}, \vz^{(i)}) \geq \frac{\sigma}{2}\|\nabla \phi^*(\vz^{(i-1)}) - \nabla\phi^*(\vz^{(i)})\|^2$. Therefore, for the quadratic term in (\ref{eq:amd-f-particular-x}) to cancel the remaining term in  (\ref{eq:amd-discr-err}), $D_{\phi^*}(\vz^{(i-1)}, \vz^{(i)})$, it suffices to have $\frac{{a_i}^2}{A^{(i)}}\leq \frac{\sigma}{L}$. It is easy to verify that $a_i = \frac{\sigma}{L}\frac{i+1}{2}$ from the theorem statement satisfies $\frac{{a_i}^2}{A^{(i)}}\leq \frac{\sigma}{L}$, and thus it follows that $G^{(k)}\leq \frac{a_0G^{(0)}}{A^{(k)}}$.

It remains to bound the initial gap, while the final bound will follow by simple computation of $A^{(k)}$. We have:
\begin{align*}
a_0G^{(0)} =& a_0(f(\vxh^{(0)})-f(\vx^{(0)})) - a_0\innp{\nabla f(\vx^{(0)}), \nabla \phi^*(\vz^{(0)})-\vx^{(0)}}\\
&- D_{\psi}(\nabla\phi^*(\vz^{(0)}), \vx^{(0)}) + D_{\psi}(\vx^*, \vx^{(0)})\\
\leq & D_{\psi}(\vx^*, \vx^{(0)}),
\end{align*}
by the same arguments as in bounding the discretization error above.
\end{proof}
\subsection{Gradient Descent} The discretization error of gradient descent is the same as the discretization error of the accelerated method from the previous subsection (Eq.~(\ref{eq:amd-fe-discr-error})), since the two methods use the same approximate duality gap. Classical gradient descent uses Forward Euler discretization, which sets $\vx^{(i)} = \vx^{(0)} - \nicefrac{\vz^{(i-1)}}{\sigma}$. Thus, the algorithm can be stated as:
\begin{equation}\label{eq:grad-desc}\tag{GD}
\begin{gathered}
\vz^{(i)} = \vz^{(i-1)} - a_i \nabla f(\vx^{(i)}),\\
\vx^{(i)} = \nabla \phi^*(\vz^{(i-1)}) = \vx^{(0)} + \frac{\vz^{(i-1)}}{\sigma},\\
\vz^{(0)} = 0, \; \vx^{(0)} \in \mathbb{R}^n \text{ is an arbitrary initial point.}
\end{gathered}
\end{equation}

To cancel out the discretization error, we only need to use the fact that gradient steps (corresponding to the steps of the algorithm) reduce the function value, assuming that the function is $L$-smooth for some $L \in \mathbb{R}_{++}$. This is achieved by setting $U^{(i)} = f(\vx^{(i+1)})$. Correcting the discretization error by the change in the upper bound:
\begin{proposition}\label{prop:gd-discr-err}
The discretization error of (\ref{eq:grad-desc}) is:
\begin{equation*}
\begin{aligned}
E_d^{(i)} =& A^{(i)}\left(f(\vx^{(i+1)}) - f(\vx^{(i)})\right) - a_i \innp{\nabla f(\vx^{(i)}), \vx^{(i+1)} - \vx^{(i)}} - \frac{{a_i}^2}{2\sigma}\|\nabla f(\vx^{(i)})\|_2^2. 
\end{aligned}
\end{equation*}
\end{proposition}
\begin{proof}
Follows directly by combining (\ref{eq:amd-fe-discr-error}), $U^{(i)} = f(\vx^{(i+1)})$, and (\ref{eq:grad-desc}).
\end{proof}
We can now recover the classical bound for gradient descent (see, e.g., \cite{Bube2014}).
\begin{theorem}\label{thm:grad-desc}
Let $f: \mathbb{R}^n \rightarrow \mathbb{R}$ be an $L$-smooth function for $L \in \mathbb{R}_{++}$, $\phi(\vx) = \frac{\sigma}{2}\|\vx - \vx^{(0)}\|^2$, and let $\vx^{(i)}$ evolve according to (\ref{eq:grad-desc}). If $a_i = \frac{\sigma}{L}$, $\forall i \geq 0$, then, $\forall k \geq 0$:
$$
f(\vx^{(k+1)}) - f(\vx^*) \leq \frac{L\|\vx^* - \vx^{(0)}\|^2}{2(k+1)}.
$$
\end{theorem}
\begin{proof}
By the smoothness of $f(\cdot)$ and (\ref{eq:grad-desc}):
\begin{align*}
E_d^{(i)} &\leq A^{(i-1)}\innp{\nabla f(\vx^{(i)}), \vx^{(i+1)}-\vx^{(i)}} + A^{(i)}\frac{L}{2}\|\vx^{(i+1)} - \vx^{(i)}\|_2^2 - \frac{{a_i}^2}{2\sigma}\|\nabla f(\vx^{(i)})\|_2^2\\
&= \left(-\frac{a_i A^{(i-1)}}{\sigma} + \frac{A^{(i)}{a_i}^2 L}{2\sigma^2} - \frac{{a_i}^2}{2\sigma}\right)\|\nabla f(\vx^{(i)})\|_2^2. 
\end{align*}
By type-checking the last expression, it follows that $a_i$ needs to be proportional to $\frac{\sigma}{L}$ to obtain $E_d^{(i)}\leq 0$. Choose $a_i = \frac{\sigma}{L}$. It remains to bound the initial gap. Observing that for $a_i = \frac{\sigma}{L}$, $U^{(0)} = f(\vx^{(0)})- \frac{1}{2L}\|\nabla f(\vx^{(0)})\|_2^2$ and $L^{(0)} = f(\vx^{(0)}) - \frac{1}{2L}\|\nabla f(\vx^{(0)})\|_2^2 + \frac{L}{2a_0}\|\vx^* - \vx^{(0)}\|_2^2$, the claimed bound on the gap follows.   
\end{proof}
\subsection{Accelerated Smooth and Strongly Convex Minimization}

Recall the accelerated dynamics for $\sigma$-strongly convex objectives (\ref{eq:ct-amd-sc}). The dynamics is almost the same as (\ref{eq:ct-amd}), except that instead of a fixed $\phi_i(\cdot)$, we now have: $\phi_i(\vx) = \int_0^i \frac{\sigma}{2}\|\vx-\vx^{(\tau)}\|^2\mathrm{d}\alpha^{(\tau)} + \phi(\vx)$. Observe that for $i \geq j$, $\phi_i(\vx)\geq \phi_j(\vx)$, $\forall \vx \in X$. We can compute the  discretization error for (\ref{eq:ct-amd-sc}) as follows. 

\begin{proposition}\label{prop:amd-sc-fe-discr-error}
The discretization error for (\ref{eq:ct-amd-sc}) is:
\begin{align*}
E_d^{(i)} \leq & A^{(i-1)}(f(\vx^{(i)})-f(\vx^{(i-1)})) - \innp{\vz^{(i)} - \vz^{(i-1)}, \vx^{(i)} - \nabla\phi_{i-1}^*(\vz^{(i)})}\notag\\
&- D_{\phi_{i-1}^*}(\vz^{(i-1)}, \vz^{(i)}). 
\end{align*}
\end{proposition}
\begin{proof}
To compute the discretization error, we first need to compute $I_c^{(i-1, i)}$:
\begin{align}
I_c^{(i-1, i)} =& -\int_{i-1}^i \innp{\nabla f(\vx^{(\tau)}), \nabla \phi_i^*(\vz^{(i)})-\vx^{(\tau)}}\mathrm{d}\alpha^{(\tau)}\notag\\
=& -\int_{i-1}^i \innp{\nabla f(\vx^{(\tau)}), \alpha^{(\tau)}\dot{\vx}^{(\tau)}}d\tau + \int_{i-1}^i \innp{\dot{\vz}^{(\tau)}, \nabla \phi_i^*(\vz^{(i)})}d\tau\notag\\
&- \int_{i-1}^i \innp{\dot{\vz}^{(\tau)}, \nabla \phi_{\tau}^*(\vz^{(\tau)})}d\tau\notag\\
=& -A^{(i-1)}(f(\vx^{(i)})-f(\vx^{(i-1)})) + \innp{\vz^{(i)}-\vz^{(i-1)}, \nabla\phi_i^*(\vz^{(i)})}\label{eq:amd-sc-I-c}\\
&- \phi_i^*(\vz^{(i)})+\phi_{i-1}^*(\vz^{(i-1)}).\notag
\end{align}
Combining (\ref{eq:I-i-1-i}), (\ref{eq:amd-sc-I-c}), and the fact that $-a_i\nabla f(\vx^{(i)}) = \vz^{(i)}-\vz^{(i-1)}$:
\begin{equation}
E_d^{(i)} = A^{(i-1)}(f(\vx^{(i)})-f(\vx^{(i-1)})) - \innp{\vz^{(i)} - \vz^{(i-1)}, \vx^{(i)}} + \phi_i^*(\vz^{(i)}) - \phi_{i-1}^*(\vz^{(i-1)}). \label{eq:amd-sc-discr-err-gen}
\end{equation}
As $\phi_i(\vx)\geq \phi_{i-1}(\vx)$, $\forall \vx \in X$, it follows that also $\phi_i^*(\vz) \leq \phi_{i-1}^*(\vz)$, $\forall \vz$. Using the definition of Bregman divergence and convexity of $f(\cdot)$:
\begin{align}
E_d^{(i)} \leq & A^{(i-1)}(f(\vx^{(i)})-f(\vx^{(i-1)})) - \innp{\vz^{(i)} - \vz^{(i-1)}, \vx^{(i)} - \nabla\phi_{i-1}^*(\vz^{(i)})}\notag\\
&- D_{\phi_{i-1}^*}(\vz^{(i-1)}, \vz^{(i)}).
\end{align}
\end{proof}

Comparing the discretization error from Proposition \ref{prop:amd-sc-fe-discr-error} with the discretization error (\ref{eq:amd-fe-discr-error}) from previous subsection, we can observe that they take the same form, with the only difference of $\phi^*$ being replaced by $\phi_{i-1}^*$. Thus, introducing a gradient descent step into the discrete algorithm leads to the same changes in the discretization error, and we can use the same arguments to analyze the convergence. The algorithm is given as:
\begin{equation}\label{eq:amd-sc}\tag{ASC}
\begin{gathered}
\vz^{(i)} = \vz^{(i-1)} - a_i \nabla f(\vx^{(i)}),\\
\vx^{(i)} = \frac{A^{(i-1)}}{A^{(i)}}\vxh^{(i-1)} + \frac{a_i}{A^{(i)}}\nabla\phi_{i-1}^*(\vz^{(i-1)}),\\
\vxh^{(i)} = \mathrm{Grad}(\vx^{(i)}),\\
\vz^{(0)} = -a_0 f(\vx^{(0)}),\; \vxh^{(0)} = \mathrm{Grad}(\vx^{(0)}), \text{ for arbitrary }\vx^{(0)} \in X,
\end{gathered}
\end{equation}
while the discretization error for (\ref{eq:amd-sc}) becomes:
\begin{equation}\label{eq:amd-sc-discr-err}
\begin{aligned}
E_d^{(i)} &\leq  A^{(i)}(f(\vxh^{(i)}) - f(\vx^{(i)})) + a_i \innp{\nabla f(\vx^{(i)}), \nabla \phi_{i-1}^*(\vz^{(i-1)}) - \nabla \phi_{i-1}^*(\vz^{(i)})}\\
&- D_{\phi_{i-1}^*}(\vz^{(i-1)}, \vz^{(i)}). 
\end{aligned}
\end{equation}
We have the following convergence result:
\begin{theorem}\label{thm:amd-sc-conv}
Let $f:X\rightarrow \mathbb{R}$ be an $L$-smooth and $\sigma$-strongly convex function, $\psi:X\rightarrow\mathbb{R}$ be a $\sigma_0$-strongly convex function, for $\sigma_0 = L - \sigma$, $\phi(\cdot) = D_{\psi}(\cdot, \vx^{(0)})$, and let $\vx^{(i)}, \vxh^{(i)}, \vz^{(i)}$ evolve according to (\ref{eq:amd-sc}), where $\phi_i(\vx) = \sum_{j=0}^i a_j\frac{\sigma}{2}\|\vx - \vx^{(j)}\|^2 + \phi(\vx)$, $\forall \vx \in X$. If $a_0 = 1$ and $\frac{a_i}{A^{(i)}} = \frac{\sqrt{4\kappa + 1} - 1}{2\kappa}$, where $\kappa = L/\sigma$ is $f(\cdot)$'s condition number, then:
$$
f(\vxh^{(k)}) - f(\vx^*) \leq \left(1-\frac{\sqrt{4\kappa + 1} - 1}{2\kappa}\right)^k D_{\psi}(\vx^*, \vx^{(0)}).
$$
\end{theorem}
\begin{proof}
The proof follows by applying the same arguments as in the proof of Theorem \ref{thm:amd-conv}. To obtain the convergence bound, we observe that $\phi_i(\cdot)$ is $\sigma_i$-strongly convex for $\sigma_i = \sigma\sum_{j=0}^ia_j + \sigma_0 = A^{(i)}\sigma + \sigma_0$. Thus, we only need to show that $\frac{{a_i}^2}{A^{(i)}}\leq \frac{\sigma_{i-1}}{L}$. A sufficient condition is that $\frac{{a_i}^2}{A^{(i)}A^{(i-1)}}\leq \frac{\sigma}{L} = \frac{1}{\kappa}$, which is equivalent to $\frac{{a_i}^2}{(A^{(i)})^2}\leq \frac{1}{\kappa}(1-\frac{a_i}{A^{(i)}})$. Solving $\frac{{a_i}^2}{(A^{(i)})^2} = \frac{1}{\kappa}(1-\frac{a_i}{A^{(i)}})$ gives the $a_i$'s from the theorem statement for $i\geq 1$. The choice of $a_0 = 1$, $\sigma_0 = L-\sigma$ ensures $a_0G^{(0)} \leq \phi(\vx^*) = D_{\psi}(\vx^*, \vx^{(0)})$.
\end{proof}

\begin{remark}
When $X = \mathbb{R}^n$, we obtain a tighter convergence bound. Namely, assuming that $\|\cdot\| = \|\cdot\|_2$, we can recover the standard guarantee $f(\vxh^{(k)}) - f(\vx^*) \leq \left(1-\frac{1}{\sqrt{\kappa}}\right)^k D_{\psi}(\vx^*, \vx^{(0)})$ \cite{nesterov2013introductory}. More details can be found in Appendix~\ref{sec:amd-sc-unconstrained}. 
\end{remark}

\subsection{Composite Dual Averaging}

Consider the forward Euler discretization of (\ref{eq:ct-comp-md}), recovering updates similar to \cite{duchi2010composite}\footnote{(\ref{eq:comp-md}) is the ``lazy'' (dual averaging) version of the COMID algorithm from \cite{duchi2010composite}.}:
\begin{equation}\label{eq:comp-md}\tag{CMD}
\begin{gathered}
\vz^{(i)} = \vz^{(i-1)}- a_i \nabla f(\vx^{(i)}),\\
\vx^{(i)} = \nabla\phi_{i}^*(\vz^{(i-1)}),\\
\vxh^{(i)} = \frac{A^{(i-1)}}{A^{(i)}} + \frac{a_i}{A^{(i)}}\vx^{(i)},\\
\vz^{(0)} = -a_0 f(\vx^{(0)}),\; \vxh^{(0)} = \vx^{(0)}, \text{ for arbitrary } \vx^{(0)} \in X. 
\end{gathered}
\end{equation}

Unlike in the standard (non-composite) convex minimization, as discussed at the beginning of the section, in the composite case the discretization error needs to take into account an extra term. The additional term appears due to the discontinuous solution updates and $\psi(\cdot)$ in the objective; in the continuous-time case $\vx^{(t)}=\nabla \phi_t^*(\vz^{(t)})$ and the change in the upper bound term $\int_{0}^t\psi(\vx^{(\tau)})\mathrm{d}\alpha^{(\tau)}$ matches the change in $\psi(\nabla\phi_t^*(\vz^{(t)}))$. In the discrete time, however, $\vx^{(i)} = \nabla\phi_{i}^*(\vz^{(i-1)})$, and thus the error also includes $a_i(\psi(\nabla\phi_{i}^*(\vz^{(i-1)}))-\psi(\nabla\phi_{i}^*(\vz^{(i)})))$, leading to the following bound on the discretization error.

\begin{proposition}\label{prop:comp-md-discr-err}
The discretization error for forward Euler discretization of (\ref{eq:ct-comp-md}) is:
$$
E_d^{(i)} \leq D_{\phi_i^*}(\vz^{(i)}, \vz^{(i-1)}).
$$
\end{proposition}
\begin{proof}
In continuous-time regime, $\vx^{(\tau)} = \nabla\phi_{\tau}^*(\vz^{(\tau)})$, and thus:
\begin{align*}
I_c^{(i-1, i)} =& \int_{i-1}^i \innp{\dot{\vz}^{(\tau)}, \nabla\phi_i^*(\vz^{(i)})-\nabla\phi_{\tau}^*(\vz^{(\tau)})}d\tau\\
=& \innp{\nabla\phi_i^*(\vz^{(i)}), \vz^{(i)}-\vz^{(i-1)}} - \int_{i-1}^i \innp{\nabla\phi^*_{\tau}(\vz^{(\tau)}), \dot{\vz}^{(\tau)}}d\tau\\
=& \innp{\nabla\phi_i^*(\vz^{(i)}), \vz^{(i)}-\vz^{(i-1)}} - \int_{i-1}^i \left(\frac{\d}{\d\tau}\phi_{\tau}^*(\vz^{(\tau)}) - \frac{d}{d s}\left.\phi_{s}^*(\vz^{(\tau)})\right|_{s=\tau}\right)d\tau\\
=& \innp{\nabla\phi_i^*(\vz^{(i)}), \vz^{(i)}-\vz^{(i-1)}} - \phi_i^*(\vz^{(i)}) + \phi_{i-1}^*(\vz^{(i-1)})\\
&+ \int_{i-1}^{i}\frac{d}{d s}\left.\phi_{s}^*(\vz^{(\tau)})\right|_{s=\tau}d\tau.
\end{align*}
Recalling that $\phi_t(\vx) = \alpha^{(t)}\psi(\vx) + \phi(\vx)$ and using Danskin's theorem, we have:
\begin{align*}
\int_{i-1}^{i}\frac{d}{d s}\left.\phi_{s}^*(\vz^{(\tau)})\right|_{s=\tau}d\tau &= \int_{i-1}^i \frac{d}{d s}\left.\max_{\vx \in X}\left\{\innp{\vz^{(\tau)}, \vx}-\alpha^{(s)}\psi(\vx)-\phi(\vx)\right\}\right|_{s=\tau}d\tau\\
&= -\int_{i-1}^i \dot{\alpha}^{(\tau)}\psi(\nabla\phi^*_{\tau}(\vz^{(\tau)}))d\tau = -a_i\psi(\nabla\phi^*_i(\vz^{(i)})).
\end{align*}
Therefore:
$$
I_c^{(i-1, i)} = \innp{\nabla\phi_i^*(\vz^{(i)}), \vz^{(i)}-\vz^{(i-1)}} - \phi_i^*(\vz^{(i)}) + \phi_{i-1}^*(\vz^{(i-1)})-a_i\psi(\nabla\phi^*_i(\vz^{(i)})).
$$
On the other hand, as:
\begin{align*}
I^{(i-1, i)} &= -a_i\innp{\nabla f(\vx^{(i)}), \nabla\phi_i^*(\vz^{(i)})-\vx^{(i)}}\\
&= \innp{\vz^{(i)}-\vz^{(i-1)}, \nabla\phi_i^*(\vz^{(i)}) - \nabla\phi_{i}^*(\vz^{(i-1)})}, 
\end{align*}
the discretization error is:
\begin{align*}
E_d^{(i)} &= \phi_i^*(\vz^{(i)}) - \phi_{i-1}^*(\vz^{(i-1)}) -  \innp{\nabla\phi_{i}^*(\vz^{(i-1)}), \vz^{(i)}-\vz^{(i-1)}}  + a_i\psi(\nabla\phi^*_{i}(\vz^{(i-1)}))\\
&= D_{\phi_i^*}(\vz^{(i)}, \vz^{(i-1)}) + \phi_{i}^*(\vz^{(i-1)}) - \phi_{i-1}^*(\vz^{(i-1)})+ a_i\psi(\nabla\phi^*_{i}(\vz^{(i-1)})).
\end{align*}
It remains to show that $\phi_{i}^*(\vz^{(i-1)}) - \phi_{i-1}^*(\vz^{(i-1)})+ a_i\psi(\nabla\phi^*_{i}(\vz^{(i-1)}))\leq 0$. Observing that $a_i\psi(\nabla\phi^*_{i}(\vz^{(i-1)})) = \phi_i(\nabla \phi_i^*(\vz^{(i-1)}))-\phi_{i-1}(\nabla \phi_i^*(\vz^{(i-1)}))$ and using the definition of a convex conjugate together with Fact \ref{fact:danskin}:
\begin{align*}
&\phi_{i}^*(\vz^{(i-1)}) - \phi_{i-1}^*(\vz^{(i-1)})+ a_i\psi(\nabla\phi^*_{i}(\vz^{(i-1)}))\\
&\hspace{1cm}= \phi_{i}^*(\vz^{(i-1)}) - \phi_{i-1}^*(\vz^{(i-1)}) + \phi_i(\nabla \phi_i^*(\vz^{(i-1)}))-\phi_{i-1}(\nabla \phi_i^*(\vz^{(i-1)}))\\
&\hspace{1cm}= \innp{\vz^{(i-1)}, \nabla\phi_i^*(\vz^{(i-1)})-\nabla\phi_{i-1}^*(\vz^{(i-1)})}\\
&\hspace{1.3cm}+\phi_{i-1}(\nabla \phi_{i-1}^*(\vz^{(i-1)}))-\phi_{i-1}(\nabla \phi_{i}^*(\vz^{(i-1)}))\\
&\hspace{1cm}\leq 0,
\end{align*}
where the inequality is by Fact \ref{fact:danskin}, as $\nabla\phi_{i-1}(\vz^{(i-1)})=\arg\min_{\vx\in X}\{-\innp{\vz^{(i-1)}, \vx} + \phi_{i-1}(\vx)\}.$
\end{proof}

Finally, we can obtain the following convergence result for the composite functions, similar to the classical case of mirror descent.
\begin{theorem}\label{thm:composite-md}
Let $\bar{f} = f + \psi: X\rightarrow \mathbb{R}$ be a composite function, such that $f(\cdot)$ is $L$-Lipschitz-continuous and convex, and $\psi(\cdot)$ is ``simple'' and convex. Here, ``simple'' means that $\nabla\phi_i^*(\vz)$ is easily computable for $\phi_i(\cdot) = A^{(i)}\psi(\cdot) + D_\phi(\cdot, \vx^{(0)})$ and some $\sigma$-strongly convex $\phi(\cdot)$ where $\sigma > 0$. Fix any $k\geq 1$ and let $\vx^{(i)}, \vxh^{(i)}$ evolve according to (\ref{eq:comp-md}) for $a_i = \frac{1}{L}\sqrt{\frac{2\sigma \phi(\vx^*)}{k+1}}$. Then:
$$
\bar{f}(\vxh^{(k)}) - \bar{f}(\vx^*) \leq \sqrt{\frac{2D_{\phi}(\vx^*, \vx^{(0)})}{\sigma}}\frac{L}{\sqrt{k+1}}.
$$
\end{theorem}
\begin{proof}
Observe that since $\phi(\cdot)$ is $\sigma$-strongly convex, $\phi_i(\cdot)$ is also $\sigma$-strongly convex. The rest of the proof follows the same argument as the proof of Theorem \ref{thm:md} (dual averaging convergence), as $D_{\phi_i^*}(\vz^{(i)}, \vz^{(i-1)}) = -a_i\innp{\nabla f(\vx^{(i)}), \vx^{(i+1)}-\vx^{(i)}}-D_{\phi_i^*}(\vz^{(i-1)}, \vz^{(i)})$, and is omitted.
\end{proof}
\begin{remark}
Observe that if $\psi(\cdot)$ was $\sigma$-strongly convex for some $\sigma > 0$, we could have obtained a stronger convergence result, as in that case we would have $D_{\phi_i^*}(\vz^{(i-1)}, \vz^{(i)})\geq A^{(i)}\frac{\sigma}{2}\|\vx^{(i)}-\vx^{(i+1)}\|^2$, which would allow choosing larger steps $a_i$.
\end{remark}

\subsection{Frank-Wolfe Method}

For the discretization of continuous-time Frank-Wolfe method (\ref{eq:ct-fw}), we need to take into account the different lower bound we obtained in Equation~\eqref{eq:gen-lb-fw}. In particular, the integral that accrues a discretization error is $-\int_{i-1}^i \innp{\nabla f(\vx^{(\tau)}), \nabla\psi^*(\vzh^{(\tau)})-\vx^{(\tau)}}\mathrm{d}\alpha^{(\tau)}$, where $\vzh^{(\tau)} = -\nabla f(\vx^{(\tau)})$. The forward Euler discretization gives the following algorithm:
\begin{equation}\label{eq:fw}\tag{FW}
\begin{gathered}
\vzh^{(i)} = - \nabla f(\vx^{(i)}),\\
\vx^{(i)} = \frac{A^{(i-1)}}{A^{(i)}}\vx^{(i-1)} + \frac{a_i}{A^{(i)}}\nabla\psi^*(\vzh^{(i-1)}),\\
\vx^{(0)} \in X \text{ is an arbitrary initial point}.
\end{gathered}
\end{equation} 
As discussed before, the discretization error needs to include  $a_i(\psi(\nabla\psi^*(\vzh^{(i-1)}))-\psi(\nabla\psi^*(\vzh^{(i)})))$ in addition to $I^{(i-1, i)}-I_c^{(i-1, i)}$, and is bounded as follows. 
\begin{proposition}\label{prop:fw-discr-err}
The discretization error for forward Euler discretization of (\ref{eq:ct-fw}) is:
$$
E_d^{(i)} \leq a_i \innp{\nabla f(\vx^{(i)})-\nabla f(\vx^{(i-1)}), \nabla\psi^*(\vzh^{(i-1)})-\nabla\psi^*(\vzh^{(i)})}. 
$$
\end{proposition}
\begin{proof}
In the discrete-time case, $I^{(i-1, i)} = -a_i\innp{\nabla f(\vx^{(i)}), \nabla\psi^*(\vzh^{(i)})-\vx^{(i)}}$, while in continuous-time case, as $\dot{\vx}^{(t)} = \dot{\alpha}^{(t)}\frac{\nabla\psi^*(\vzh^{(t)})-\vx^{(t)}}{\alpha^{(t)}}$ and integrating by parts:
$$
I_c^{(i-1, i)} = -\int_{i-1}^{i}\innp{\nabla f(\vx^{(\tau)}), \alpha^{(\tau)}\dot{\vx}^{(\tau)}}d\tau = -A^{(i-1)}(f(\vx^{(i)})-f(\vx^{(i-1)})).
$$
Therefore:
\begin{equation}\label{eq:fw-discr-err-1}
\begin{aligned}
E_d^{(i)} =& A^{(i-1)}(f(\vx^{(i)})-f(\vx^{(i-1)})) - a_i\innp{\nabla f(\vx^{(i)}), \nabla\psi^*(\vzh^{(i)})-\vx^{(i)}}\\
&+ a_i(\psi(\nabla\psi^*(\vzh^{(i-1)}))-\psi(\nabla\psi^*(\vzh^{(i)})))\\
\leq & a_i\innp{\nabla f(\vx^{(i)}), \nabla\psi^*(\vzh^{(i-1)})-\nabla\psi^*(\vzh^{(i)})}\\
&+ a_i(\psi(\nabla\psi^*(\vzh^{(i-1)}))-\psi(\nabla\psi^*(\vzh^{(i)}))), 
\end{aligned}
\end{equation}
where we have used the convexity of $f(\cdot)$ and $\vx^{(i)}=\frac{A^{(i-1)}}{A^{(i)}}\vx^{(i-1)}+\frac{a_i}{A^{(i)}}\nabla\psi^*(\vzh^{(i-1)})$. 
Further, by Fact \ref{fact:danskin}, 
\begin{align}
-\innp{\vzh^{(i-1)}, \nabla\psi^*(\vzh^{(i-1)})} + \psi(\nabla\psi^*(\vzh^{(i-1)}))\leq -\innp{\vzh^{(i-1)}, \nabla\psi^*(\vzh^{(i)})} + \psi(\nabla\psi^*(\vzh^{(i)})),\notag
\end{align}
and, therefore, as $\vzh^{(i-1)} = -\nabla f(\vx^{(i-1)})$:
\begin{align}
\psi(\nabla\psi^*(\vzh^{(i-1)})) - \psi(\nabla\psi^*(\vzh^{(i)})) \leq \innp{\nabla f(\vx^{(i-1)}), \nabla\psi^*(\vzh^{(i)}) - \nabla\psi^*(\vzh^{(i-1)})}. \label{eq:fw-discr-err-psi-terms}
\end{align}

Combining (\ref{eq:fw-discr-err-1}) and (\ref{eq:fw-discr-err-psi-terms}), the claimed bound on discretization error follows.
\end{proof}

We can now recover the convergence result from \cite{nesterov2015cgm}.

\begin{theorem}\label{thm:fw}
Let $\bar{f} = f+\psi: X\rightarrow \mathbb{R}$ be a composite function, where $\psi(\cdot)$ is convex and $f(\cdot)$ is convex with H\"{o}lder-continuous gradients, i.e., for some fixed $L_{\nu}<\infty$, $\nu \in (0, 1]$\footnote{Observe that when $\nu = 1$, $f(\cdot)$ is $L_{\nu}$-smooth.}:
\begin{equation}
\|\nabla f(\vx) - \nabla f(\vxh)\|_*\leq L_{\nu}\|\vx - \vxh\|^{\nu}, \; \forall \vx, \vxh \in X.\label{eq:f-holder-cont}
\end{equation}
Let $D \defeq \max_{\vx, \vxh \in X}\|\vx - \vxh\|$ denote the diameter of $X$. 
If $\vx^{(i)}$ evolves according to (\ref{eq:fw}), then, $\forall k \geq 1$:
$$
\bar{f}(\vx^{(k)}) - \bar{f}(\vx^*) \leq L_{\nu}D^{1+\nu}\frac{1}{A^{(k)}}{\sum_{i=0}^k \frac{{a_i}^{1+\nu}}{(A^{(i)})^{\nu}}}.
$$
In particular, if $a_i = i+1$, then
$$
\bar{f}(\vx^{(k)}) - \bar{f}(\vx^*) \leq 2^{1+ \nu}\frac{L_{\nu}D^{1+\nu}}{(k+1)^{\nu}}.
$$
\end{theorem}
\begin{proof}
Applying Cauchy-Schwartz Inequality to the discretization error given by Proposition \ref{prop:fw-discr-err}, we have:
\begin{align*}
E_d^{(i)} &\leq a_i \|\nabla f(\vx^{(i)}) - \nabla f(\vx^{(i-1)})\|_*\cdot \|\nabla \psi^*(\vzh^{(i)}) - \nabla \psi^*(\vzh^{(i-1)})\|\\
&\leq  \frac{{a_i}^{1+\nu}}{(A^{(i)})^{\nu}}L_{\nu}\|\vx^{(i-1)} - \nabla \psi^*(\vzh^{(i-1)})\|^{\nu}\cdot \|\nabla \psi^*(\vzh^{(i)}) - \nabla \psi^*(\vzh^{(i-1)})\|\\
&\leq \frac{{a_i}^{1+\nu}}{(A^{(i)})^{\nu}}L_{\nu} D^{1+\nu},
\end{align*}
where the second inequality follows from (\ref{eq:f-holder-cont}) and $\vx^{(i)}-\vx^{(i-1)} = \frac{a_i}{A^{(i)}}(\vx^{(i-1)} - \nabla \psi^*(\vzh^{(i-1)}))$ (by (\ref{eq:fw})). Therefore: $G^{(k)}\leq \frac{a_0 G^{(0)}}{A^{(k)}} + L_{\nu}D^{1+\nu}\frac{1}{A^{(k)}}{\sum_{i=1}^k \frac{{a_i}^{1+\nu}}{(A^{(i)})^{\nu}}}$. 

We now use the same arguments to bound $G^{(0)}$. As $\vx^{(0)}$ can be mapped to $\nabla\psi^*(\vzh^{(-1)})$, for some $\vzh^{(-1)}$, we have:
\begin{align*}
G^{(0)} &= - \innp{\nabla f(\vx^{(0)}), \nabla \psi^*(\vzh^{(0)})-\vx^{(0)}} - \psi(\nabla\psi^*(\vzh^{(0)})) + \psi(\vx^{(0)})\leq {L_{\nu}}D^{1+\nu}.
\end{align*}
Therefore, $G^{(k)}\leq L_{\nu}D^{1+\nu}\frac{1}{A^{(k)}}{\sum_{i=0}^k \frac{{a_i}^{1+\nu}}{(A^{(i)})^{\nu}}}$. In particular, when $a_i = i+1$, then $A^{(i)} = \frac{(i+1)(i+2)}{2}$. Finally, $\sum_{i=0}^k \frac{{a_i}^{1+\nu}}{(A^{(i)})^{\nu}} < 2^{\nu}\sum_{i=0}^k (i+1)^{1-\nu}< 2^{\nu}(k+1)^{2-\nu}$, and the convergence bound follows. 
\end{proof}

\section{Conclusion}\label{sec:conclusion}
We presented a general technique for the analysis of first-order methods. The technique is intuitive and follows the argument of reducing approximate duality gap at the rate equal to the convergence rate. Besides the unified interpretation of many first-order methods, the technique is generally useful for obtaining new optimization methods~\cite{AXGD,LP-jelena-lorenzo,diakonikolas2018fairpc,cohen2018acceleration,diakonikolas2018alternating}. An interesting direction for future is extending this technique to other settings, such as, e.g., geodesically convex optimization.

\section*{Acknowledgements}

We thank the anonymous reviewers for their thoughtful comments and suggestions, which have greatly improved the presentation of this paper. We also thank Ziye Tang for pointing out several typos in the earlier version of the paper and providing useful suggestions for improving its presentation. 

\bibliographystyle{siamplain}
{\small
\bibliography{references}
}
\appendix
\section{Properties of the Bregman  Divergence}\label{sec:breg-div-prop}

The following properties of Bregman divergence are useful in our analysis.

\begin{proposition}\label{prop:cvx-conj-bd-is-strongly-cvx-too}
If $\psi(\cdot)$ is $\sigma$-strongly convex, then $D_{\psi^*}(\vz, \vzh) \geq \frac{\sigma}{2}\|\nabla \psi^*(\vz) - \nabla \psi^*(\vzh)\|^2$.
\end{proposition}
\begin{proof}
From the definition of $\psi^*$ and Fact \ref{fact:danskin}, 
\begin{equation}\label{eq:psi*-in-terms-of-psi}
 \psi^*(\vz) = \innp{\nabla \psi^*(\vz), \vz} - \psi(\nabla \psi^*),\; \forall \vz.
 \end{equation}
Using the definition of $D_{\psi^*}(\vz, \vzh)$ and (\ref{eq:psi*-in-terms-of-psi}), we can write $D_{\psi^*}(\vz, \vzh)$ as:
$$
D_{\psi^*}(\vz, \vzh) = \psi(\nabla \psi^*(\vzh)) - \psi(\nabla\psi^*(\vz)) - \innp{\vz, \nabla \psi^*(\vzh) - \nabla \psi^*(\vz)}.
$$
Since $\psi(\cdot)$ is $\sigma$-strongly convex, it follows that:
$$
D_{\psi^*}(\vz, \vzh) \geq \frac{\sigma}{2}\|\nabla \psi^*(\vzh) - \nabla \psi^*(\vz)\|^2 + \innp{\nabla \psi(\nabla \psi^*(\vz))-\vz, \nabla \psi^*(\vzh) - \nabla \psi^*(\vz)}.
$$
As, from Fact \ref{fact:danskin}, $\nabla\psi^*(\vz) = \arg\max_{\vx \in X}\{\innp{\vx, \vz} - \psi(\vx)\}$, by the first-order optimality condition $$\innp{\nabla \psi(\nabla \psi^*(\vz))-\vz, \nabla \psi^*(\vzh) - \nabla \psi^*(\vz)} \geq 0,$$ completing the proof.
\end{proof}

The Bregman divergence $D_{\psi^*}(\vx,\vy)$ captures the difference between $\psi^*(\vx)$ and its first order approximation at $\vy.$ Notice that, for a differentiable $\psi^*$, we have:
$$
\nabla_{\vx} D_{\psi^*}(\vx,\vy) = \nabla \psi^*(\vx) - \nabla \psi^*(\vy).
$$
The Bregman divergence $D_{\psi^*}(\vx,\vy)$ is a convex function of $\vx.$ Its Bregman divergence is itself.
\begin{proposition}\label{prop:magic-identity}
For all $\vx, \vy, \vz \in X$ 
$$
D_{\psi^*}(\vx, \vy) = D_{\psi^*}(\vz, \vy) + \innp{\nabla \psi^*	(\vz) - \nabla \psi^*(\vy), \vx - \vz} + D_{\psi^*}(\vx, \vz).
$$
\end{proposition}
\section{Extension of ADGT to Monotone Operators and Convex-Concave Saddle-Point Problems}\label{app:mon-ops}
Given a monotone operator $F:X\rightarrow \mathbb{R}^n$, the goal is to find a point $\vx^*\in X$ such that $\innp{F(\vu), \vx^* - \vu} \leq 0$, $\forall \vu \in X$. The approximate version of this problem is:
\begin{equation}\label{eq:approx-mon-op}
\text{Find } \vx_{\epsilon}\in X \text{ such that } \innp{F(\vu), \vx_{\epsilon} - \vu} \leq \epsilon, \quad \forall \vu \in X,
\end{equation}
and we can think of $\epsilon$ on the right-hand side of (\ref{eq:approx-mon-op}) as the optimality gap.

The property of monotone operators $F(\cdot)$ useful for the approximate gap analysis is $\forall \vx, \vu \in X$: $\innp{F(\vu), \vx - \vu}\leq \innp{F(\vx), \vx - \vu}$.  
The approximate gap can be constructed using the same ideas as in the case of a convex function, which, letting $\vxh^{(t)} = \frac{1}{\alpha^{(t)}}\int_{0}^t \vx^{(\tau)}\mathrm{d}\alpha + \frac{\alpha^{(t)}-A^{(t)}}{\alpha^{(t)}}\vx^{(0)}$,  
gives, $\forall \vu \in X$:
\begin{equation}\label{eq:mon-op-gap}
\begin{aligned}
&\innp{F(\vu), \vxh^{(t)} - \vu} \leq G^{(t)}\\
&\hspace{1cm}\defeq \frac{\max_{\vx\in X}\left\{ \int_{0}^t \innp{F(\vx^{(\tau)}), \vx^{(\tau)} - \vx}\mathrm{d}\alpha - \phi(\vx) \right\} + \max_{\vu \in X}\phi(\vu)}{\alpha^{(t)}}\\
&\hspace{1.4cm}+  \frac{\alpha^{(t)}-A^{(t)}}{\alpha^{(t)}}\max_{\vv\in X}\innp{F(\vv), \vx^{(0)}-\vv}, 
\end{aligned}
\end{equation}

Now assume that we want to find a saddle point of a function $\Phi(\vv, \vw): V\times W \rightarrow \mathbb{R}$ that is convex in $\vv$ and concave in $\vw$. By convexity in $\vv$ and concavity in $\vw$, we have that, for all $\vv, \vv^{(\tau)} \in Y$ and all $\vw, \vw^{(\tau)} \in Z$:
\begin{align}
\Phi(\vv, \vw^{(\tau)}) - \Phi(\vv^{(\tau)}, \vw^{(\tau)}) \geq \innp{\nabla_{\vv}\Phi(\vv^{(\tau)}, \vw^{(\tau)}), \vv - \vv^{(\tau)}}, \label{eq:saddle-point-cvxity}\\
\Phi(\vv^{(\tau)}, \vw) - \Phi(\vv^{(\tau)}, \vw^{(\tau)}) \leq \innp{\nabla_{\vw}\Phi(\vv^{(\tau)}, \vw^{(\tau)}), \vw - \vw^{(\tau)}},
\label{eq:saddle-point-cncvity}
\end{align}
where $\nabla_{\vv}$ (resp. $\nabla_{\vw}$) denotes the gradient w.r.t.~$\vv$ (resp.~$\vw$). 

Combining (\ref{eq:saddle-point-cvxity}) and (\ref{eq:saddle-point-cncvity}), it follows that, $\forall \vv, \vv^{(\tau)} \in Y$, $\forall \vw, \vw^{(\tau)} \in Z$:
\begin{equation*}
\begin{aligned}
&\Phi(\vv^{(\tau)}, \vw) - \Phi(\vv, \vw^{(\tau)})\\
&\hspace{1cm}\leq \innp{\nabla_{\vv}\Phi(\vv^{(\tau)}, \vw^{(\tau)}),  \vv^{(\tau)} - \vv}
-\innp{\nabla_{\vw}\Phi(\vv^{(\tau)}, \vw^{(\tau)}), \vw^{(\tau)} - \vw }.
\end{aligned}
\end{equation*}
Let $\vx = [\vv, \vw]^T$, $F(\vx) = [\nabla_{\vv}\Phi(\vv, \vw), -\nabla_{\vw}\Phi(\vv, \vw)]^T$, and  $\overline{{\vv}} = \frac{1}{A^{(t)}}\int_{0}^t \vv^{(\tau)} \mathrm{d}\alpha^{(\tau)}$, $\overline{\vw} = \frac{1}{A^{(t)}}\int_{0}^t \vw^{(\tau)}\mathrm{d}\alpha^{(\tau)}$. Then, we have, $\forall \vx = [\vv, \vw]^T\in V\times W$:
\begin{equation*}
\Phi(\overline{\vv}, \vw) - \Phi(\vv, \overline{\vw}) \leq \frac{\int_{0}^t \innp{F(\vx^{(\tau)}), \vx^{(\tau)}-\vx}\mathrm{d}\alpha^{(\tau)}}{A^{(t)}},
\end{equation*}
and using the same arguments as before, we obtain the same bound for the gap as (\ref{eq:mon-op-gap}). Therefore, we can focus on analyzing the decrease of $G^{(t)}$ from (\ref{eq:mon-op-gap}) as a function of $t$ and the same result will follow for the gap of convex-concave saddle-point problems.

\subsection{Continuous-time Mirror Descent}

Replacing $\nabla f(\vx^{(t)})$ by $F(\vx^{(t)})$ in the gap for mirror descent for convex minimization from Section~\ref{sec:ct-md}, it follows that (\ref{eq:ct-md}) also ensures $\frac{\d}{\d t}(\alpha^{(t)}G^{(t)}) = 0$ for the gap (\ref{eq:mon-op-gap}) derived for monotone operators and saddle-point problems. Hence, we have the following Lemma:
\begin{lemma}\label{lemma:ct-md-for-mon-ops}
Suppose we are given a variational inequality problem with monotone operator $F:\mathbb{R}^n\rightarrow \mathbb{R}^n$. Then a version of (\ref{eq:ct-md}) that replaces $\nabla f(\vx^{(t)})$ by $F(\vx^{(t)})$ ensures that, $\forall t > 0$, $\forall \vu \in X$:
$$
\innp{F(\vu), \vxh^{(t)}-\vu} \leq \frac{\max_{\vx' \in X}\phi(\vx') + \alpha^{(0)}\max_{\vx'' \in X}\innp{F(\vx''), \vx^{(0)}-\vx''}}{\alpha^{(t)}}.
$$
Moreover, for a convex-concave saddle-point problem $\min_{\vv \in V}\max_{\vw \in W}\Phi(\vv, \vw)$, taking $\vx = [\vv, \vw]^T$, $F(\vx) = [\nabla_{\vv} \Phi(\vv, \vw), - \nabla_{\vw}\Phi(\vv, \vw)]^T$, then the version of (\ref{eq:ct-md}) that uses the monotone operator $F(\vx)$ ensures that, $\forall t > 0$, $\forall (\vv, \vw) \in V\times W$:
\begin{align*}
 &\Phi(\vvh^{(t)}, \vw) - \Phi(\vv, \vwh^{(t)})\\
 &\hspace{1cm}\leq \frac{\max_{\vx' \in X}\phi(\vx') + \alpha^{(0)} \max_{\vv'' \in V, \vw'' \in W}\left\{\Phi(\vvh^{(0)}, \vw'') - \Phi(\vv'', \vwh^{(0)})\right\}}{\alpha^{(t)}}.
\end{align*}
\end{lemma}
\subsection{Lazy Mirror Prox}

As discussed in Section~\ref{sec:discretization}, a lazy version of the mirror prox method~\cite{Mirror-Prox-Nemirovski} follows as an approximate backward Euler (or predictor-corrector) discretization of the mirror descent dynamics~\eqref{eq:ct-md}, stated in~\eqref{eq:mp}.  

\begin{theorem}\label{thm:mp-conv} Let $F:X\rightarrow \mathbb{R}^n$ be an $L$-smooth monotone operator and let $\psi(\cdot)$ be a $\sigma$-strongly convex function. Let $\vx^{(i)}, \vxh^{(i)}$ evolve according to (\ref{eq:mp}), where $\nabla f (\cdot)$ is replaced by $F(\cdot)$. If $a_i = \sigma/L$ and $\phi(\cdot) = D_{\psi}(\cdot, \vxt^{(0)})$, then $\forall k \geq 1$ and $\forall \vu \in X$: 
$$
\innp{F(\vu), \vxh^{(k)}-\vu} \leq \frac{L}{\sigma}\cdot\frac{\max_{\vx \in X}D_{\psi}(\vx, \vxt^{(0)})}{k}.
$$
\end{theorem}
\begin{proof}
From (\ref{eq:I-i-1-i}) and similarities between mirror-descent gaps for convex functions and monotone operators, we have that the discretization error is:
\begin{equation}\label{eq:mp-discr-error}
\begin{aligned}
E_d^{(i)} =& -a_i\innp{F(\vx^{(i)}), \nabla\phi^*(\vz^{(i)})-\vx^{(i)}} - D_{\phi^*}(\vz^{(i-1)}, \vz^{(i)})\\
=& a_i\innp{F(\vxt^{(i-1)})-F(\vx^{(i)}), \nabla\phi^*(\vz^{(i)})-\vx^{(i)}}\\
&- a_i\innp{F(\vxt^{(i-1)}), \nabla\phi^*(\vz^{(i)})-\vx^{(i)}} - D_{\phi^*}(\vz^{(i-1)}, \vz^{(i)}).
\end{aligned}
\end{equation}
As $a_iF(\vxt^{(i-1)}) = \vz^{(i-1)}-\vzt^{(i-1)}$ and  $\vx^{(i)}=\nabla\phi^*(\vzt^{(i-1)})$, Proposition \ref{prop:magic-identity} implies:
\begin{equation}
\begin{aligned}
- a_i & \innp{F(\vxt^{(i-1)}), \nabla\phi^*(\vz^{(i)})-\vx^{(i)}} -  D_{\phi^*}(\vz^{(i-1)}, \vz^{(i)})\\
&= -D_{\phi^*}(\vz^{(i-1)}, \vzt^{(i-1)}) - D_{\phi^*}(\vzt^{(i-1)}, \vz^{(i)})\\
&\leq -\frac{\sigma}{2}\left(\|\nabla\phi^*(\vz^{(i-1)}) - \nabla\phi^*(\vzt^{(i-1)})\|^2 + \|\nabla\phi^*(\vzt^{(i-1)}) - \nabla\phi^*(\vz^{(i)})\|^2\right),
\end{aligned}
\end{equation}
where the inequality is by Proposition \ref{prop:cvx-conj-bd-is-strongly-cvx-too}.

On the other hand, by $L$-smoothness of $F(\cdot)$, using Cauchy-Schwartz Inequality:
\begin{equation}\label{eq:mp-smoothness}
\begin{aligned}
a_i & \innp{F(\vxt^{(i-1)})-F(\vx^{(i)}), \nabla\phi^*(\vz^{(i)})-\vx^{(i)}}\\
&\leq a_i L \|\vxt^{(i-1)} - \vx^{(i)}\|\cdot \|\nabla\phi^*(\vz^{(i)})-\vx^{(i)}\|\\
&= a_i L \|\nabla\phi^*(\vz^{(i-1)}) - \nabla\phi^*(\vzt^{(i-1)})\|\cdot
\|\nabla\phi^*(\vzt^{(i-1)}) - \nabla\phi^*(\vz^{(i)})\|.
\end{aligned}
\end{equation}
As $a_i = \sigma/L$, combining (\ref{eq:mp-discr-error})-(\ref{eq:mp-smoothness}) with the inequality $2ab - a^2-b^2 \leq 0$, $\forall a, b$, we get that $E_d^{(i)}\leq 0$, $\forall i$. By (\ref{eq:gap-w-discr-err}), it follows that $G^{(k)}\leq \frac{a_0 G^{(0)}}{A^{(k)}}$. 

To bound the initial gap, we use a slight modification of the gap that starts from $i=1$ instead of $i=0$. Observe that we still have $G^{(k)}\leq \frac{a_1 G^{(1)}}{A^{(k)}}$, but now $a_0 = 0$ and, therefore, $A^{(k)} = \frac{\sigma}{L}k$. As $D_{\phi}(\cdot, \cdot)=D_{\psi}(\cdot, \cdot)$:
\begin{align*}
a_1G^{(1)} &= -a_1 \innp{F(\vx^{(1)}), \nabla\phi^*(\vz^{(1)})-\vx^{(1)}} - D_{\psi}(\nabla\phi^*(\vz^{(1)}), \vxt^{(0)}) + D_{\psi}(\vx^*, \vxt^{(0)})\\
&= -a_1 \innp{F(\vx^{(1)}), \nabla\phi^*(\vz^{(1)})-\vx^{(1)}} - D_{\phi}(\nabla\phi^*(\vz^{(1)}), \vxt^{(0)}) + D_{\psi}(\vx^*, \vxt^{(0)}).
\end{align*}
Observing that $a_1 F(\vxt^{(0)}) = \vz^{(0)}-\vzt^{(0)}$, $\vx^{(1)} = \nabla\phi^*(\vzt^{(0)})$, and applying the same arguments as in bounding $E_d^{(i)}$ above, it follows that $a_1G^{(1)} \leq D_{\psi}(\vx^*, \vxt^{(0)})$.
\end{proof}
Similarly as before, as convex-concave saddle-point problems have the same gap as variational inequalities, it is straightforward to extend Theorem \ref{thm:mp-conv} to this setting (see Section \ref{sec:ct-md} and \cite{Mirror-Prox-Nemirovski}). 

\section{Unconstrained Smooth and Strongly Convex Case}\label{sec:amd-sc-unconstrained}

In this section only, we assume that $\|\cdot\| = \|\cdot\|_2$ and $\phi(\vx) = \frac{\sigma_0}{2}\|\vx - \vx^{(0)}\|^2$. From (\ref{eq:amd-sc-discr-err-gen}), regardless of the specific discretization of (\ref{eq:ct-amd-sc}), the discretization error is:
\begin{align*}
E_d^{(i)} &= A^{(i-1)}(f(\vx^{(i)})-f(\vx^{(i-1)})) + a_i\innp{\nabla f(\vx^{(i)}), \vx^{(i)}} + \phi_i^*(\vz^{(i)}) - \phi_{i-1}^*(\vz^{(i-1)}).
\end{align*}
As discussed in Section \ref{sec:amd}, incorporating a gradient step at the end of each iteration effectively adds $A^{(i)}(f(\vxh^{(i)})-f(\vx^{(i)})) - A^{(i-1)}(f(\vxh^{(i-1)})-f(\vx^{(i-1)}))$ to $E_d^{(i)}$, where $\vxh^{(i)} = \mathrm{Grad}(\vx^{(i)})$, $\forall i$, since adding a gradient step only affects the upper bound. Therefore, the discretization error for the discrete-time method obtained from (\ref{eq:ct-amd-sc}) that adds a gradient step at the end of each iteration is:
\begin{align*}
E_d^{(i)} =& A^{(i)}(f(\vxh^{(i)})-f(\vx^{(i)})) + A^{(i-1)}(f(\vx^{(i)})-f(\vxh^{(i-1)})) + a_i\innp{\nabla f(\vx^{(i)}), \vx^{(i)}}\\
&+ \phi_i^*(\vz^{(i)}) - \phi_{i-1}^*(\vz^{(i-1)})\\
\leq & A^{(i)}(f(\vxh^{(i)})-f(\vx^{(i)})) + \innp{\nabla f(\vx^{(i)}), A^{(i)} \vx^{(i)}- A^{(i-1)}\vxh^{(i-1)}}\\
&+ \phi_i^*(\vz^{(i)}) - \phi_{i-1}^*(\vz^{(i-1)}),
\end{align*}
where the inequality follows by convexity of $f(\cdot)$. 

What makes the unconstrained case special is that we can set 
\begin{equation}\label{eq:xi-amd-sc-unconstr}
\vx^{(i)} = \frac{A^{(i-1)}}{A^{(i)}}\vxh^{(i)} + \frac{a_i}{A^{(i)}}\nabla \phi_i^*(\vz^{(i-1)}).
\end{equation}
(Compare with $\vx^{(i)} = \frac{A^{(i-1)}}{A^{(i)}}\vxh^{(i)} + \frac{a_i}{A^{(i)}}\nabla \phi_{i-1}^*(\vz^{(i-1)})$ in the constrained case). This is because $\nabla \phi_i^*(\vz^{(i-1)}) = \frac{\vz^{(i-1)} + \sigma\sum_{j=0}^i a_j \vx^{(j)} + \sigma_0 \vx^{(0)}}{\sigma_i}$, where $\sigma_i = A^{(i)}\sigma + \sigma_0$. Therefore, $\vx^{(i)}$ can be determined from (\ref{eq:xi-amd-sc-unconstr}) in a closed form. This seemingly minor difference in the discretization (everything else is the same as in the case of (\ref{eq:amd-sc})) allows us to get an improved convergence bound. In particular, for $\phi(\vx) = \frac{\sigma_0}{2}\|\vx - \vx^{(0)}\|^2 = \frac{L-\sigma}{2}\|\vx - \vx^{(0)}\|^2$, we have that $\phi_i(\cdot)$ is both $\sigma_i$-smooth and $\sigma_i$-strongly convex. On the other hand, by smoothness of $f(\cdot)$, in the unconstrained case we have $f(\vxh^{(i)})-f(\vx^{(i)})\leq -\frac{1}{2L}\|\nabla f(\vx^{(i)})\|_*^2$. Therefore, we can bound the discretization error as:
\begin{align*}
E_d^{(i)} \leq & A^{(i)}(f(\vxh^{(i)})-f(\vx^{(i)})) + a_i\innp{\nabla f(\vx^{(i)}), a_i \nabla\phi_i^*(\vz^{(i-1)})}\\
&+ \phi_i^*(\vz^{(i)}) - \phi_{i-1}^*(\vz^{(i-1)})\\
\leq & -\frac{A^{(i)}}{2L}\|\nabla f(\vx^{(i)})\|_*^2 + \phi_i^*(\vz^{(i)}) - \phi_{i-1}^*(\vz^{(i-1)}) - \innp{\vz^{(i)}-\vz^{(i-1)}, \nabla\phi_i^*(\vz^{(i-1)})}\\
\leq & -\frac{A^{(i)}}{2L}\|\nabla f(\vx^{(i)})\|_*^2 + D_{\phi_i^*}(\vz^{(i)}, \vz^{(i-1)})\\
=& -\frac{A^{(i)}}{2L}\|\nabla f(\vx^{(i)})\|_*^2 + D_{\phi_i}(\nabla\phi_i^*(\vz^{(i)}), \nabla\phi_i^*(\vz^{(i-1)}))\\
\leq & -\frac{A^{(i)}}{2L}\|\nabla f(\vx^{(i)})\|_*^2 + \frac{A^{(i)}\sigma}{2}\left\| \frac{a_i\nabla f(\vx^{(i)})}{A^{(i)}\sigma} \right\|_*^2,
\end{align*}
where we have used $\phi_{i-1}^*(\vz^{(i-1)}) \geq \phi_{i}^*(\vz^{(i-1)})$ and that in the unconstrained case $D_{\phi_i^*}(\vz^{(i)}, \vz^{(i-1)}) = D_{\phi_i}(\nabla\phi_i^*(\vz^{(i)}), \nabla\phi_i^*(\vz^{(i-1)}))$. Therefore, whenever $\frac{a_i}{A^{(i)}}\leq \sqrt{\frac{\sigma}{L}}$, $E_d^{(i)}\leq 0$. Setting $\frac{a_i}{A^{(i)}} = \sqrt{\frac{\sigma}{L}} = \frac{1}{\sqrt{\kappa}}$, the improved convergence bound $f(\vxh^{(k)})-f(\vx^*)\leq \left(1 - \frac{1}{\sqrt{\kappa}}\right)^k \frac{L-\sigma}{2}\|\vx^*-\vx^{(0)}\|^2$ follows.
\end{document}

%% file: ex_shared.tex
\usepackage{lipsum}
\usepackage{amsfonts}
\usepackage{graphicx}
\usepackage{epstopdf}
\usepackage{algorithmic}
\ifpdf
  \DeclareGraphicsExtensions{.eps,.pdf,.png,.jpg}
\else
  \DeclareGraphicsExtensions{.eps}
\fi

\usepackage{dsfont}
\usepackage{cite}
\usepackage{bm}
\usepackage{xspace}
\usepackage{verbatim}
\usepackage{nicefrac}

\renewcommand{\a}{\alpha}

\renewcommand{\d}{\mathrm{d}}
\catcode`"=\active
\newcommand"[1]{\ensuremath{\sp{(#1)}}}

\DeclareMathOperator*{\argmin}{argmin}

\newcommand{\innp}[1]{\left\langle #1 \right\rangle}

\newcommand{\vx}{\mathbf{x}}
\newcommand{\vxh}{\mathbf{\hat{x}}}

\newcommand{\vzh}{\mathbf{\hat{z}}}
\newcommand{\vvh}{\mathbf{\hat{v}}}
\newcommand{\vwh}{\mathbf{\hat{w}}}

\newcommand{\vzt}{\mathbf{\tilde{z}}}
\newcommand{\vxt}{\mathbf{\tilde{x}}}
\newcommand{\vy}{\mathbf{y}}
\newcommand{\vz}{\mathbf{z}}
\newcommand{\vv}{\mathbf{v}}

\newcommand{\vu}{\mathbf{u}}
\newcommand{\vw}{\mathbf{w}}

\newcommand{\defeq}{\stackrel{\mathrm{\scriptscriptstyle def}}{=}}
\newcommand{\etal}{\textit{et al}.}

 \newsiamthm{fact}{Fact}
 \newsiamremark{remark}{Remark}
  \newsiamremark{example}{Example}

\newcommand{\axgd}{\textsc{axgd}\xspace}
\newcommand{\adgt}{\textsc{adgt}\xspace}

\numberwithin{theorem}{section}

\AtBeginDocument{%
  \paperwidth=\dimexpr
    1in + \oddsidemargin
    + \textwidth
    + 1in + \oddsidemargin
  \relax
  \paperheight=\dimexpr
    1in + \topmargin
    + \headheight + \headsep
    + \textheight
    + 1in + \topmargin
  \relax
  \usepackage[pass]{geometry}\relax
}

\newcommand{\TheTitle}{The Approximate Duality Gap Technique:\\ 
A Unified Theory of First-Order Methods} 
\newcommand{\TheAuthors}{J. Diakonikolas and L. Orecchia}

\headers{The Approximate Duality Gap Technique}{\TheAuthors}

\title{{\TheTitle}\thanks{Part of this work was done while the authors were visiting the Simons Institute for the Theory of Computing and while the first author was a postdoctoral researcher at Boston University. It was partially supported by NSF grants \#CCF-1718342 and \#CCF-1740855, by the DIMACS/Simons Collaboration on Bridging Continuous and Discrete Optimization through NSF grant \#CCF-1740425, and by DHS-ALERT subaward 505035-78050.}}

\author{
  Jelena Diakonikolas\thanks{Department of Statistics, UC Berkeley, Berkeley, CA
    (\email{jelena.d@berkeley.edu})}
  \and
  Lorenzo Orecchia\thanks{Department of Computer Science, Boston University, Boston, MA
  (\email{orecchia@bu.edu})}
}

\usepackage{amsopn}
